\documentclass{article}

\usepackage[T2A]{fontenc}
\usepackage[utf8]{inputenc}
\usepackage[english,russian]{babel}

\usepackage{amsmath}
\usepackage{amssymb}
\usepackage{mathtools}
\usepackage{graphicx}
\usepackage{xcolor}

\usepackage{csquotes}
\usepackage[backend=biber,bibencoding=utf8,sorting=none,style=gost-numeric,language=autobib,autolang=other,clearlang=true,sortcites=true,doi=false,isbn=false,date=year]{biblatex}
\usepackage[ruled]{algorithm2e}
\usepackage{algorithmic}

\newcommand{\argmin}{\operatornamewithlimits{argmin}}

\title{Алгоритмы локальной минимизации силового поля для трехмерных макромолекул}
\author{
Яковлев П.А. \\БИОКАД, Санкт-Петербург\\
\and
Аникин А.С.\\ИДСТУ СО РАН, Иркутск
\and
Большакова О.А. \\Сириус, Сочи
\and
Гасников А.В. \\БИОКАД, МФТИ, Москва
\and
Горнов А.Ю. \\ИДСТУ СО РАН, Иркутск
\and
Ермак Т.В. \\БИОКАД, Санкт-Петербург
\and
Макаренко Д.В. \\ МФТИ, Москва
\and
Морозов В.П. \\БИОКАД, Санкт-Петербург
\and
Нетеребский Б.О. \\БИОКАД, Санкт-Петербург
}
\date{\today}

\begin{document}

\maketitle

\begin{abstract}
  Большинство проблем структурной вычислительной биологии требуют решения задачи минимизации энергетической функции, определенной на геометрии молекулы.
  Это позволяет определять свойства молекул, предсказывать правильное положение белковых цепей, находить лучшую состыковку молекул при предсказании комплексообразования (докинге), проверять гипотезы относительно белкового дизайна и решать многие другие задачи, возникающие при современной разработке лекарственных средств.
  В случае низкомолекулярных соединений (состоящих из менее чем 250 атомов) задача нахождения геометрии, минимизирующей энергетическую функцию, является достаточно хорошо решенной.
  Более сложной задачей является минимизация макромолекул (в частности, белков), в состав которых входят десятки тысяч атомов.
  Однако отличительной особенностью данных постановок задач является наличие хороших начальных приближений к нужному решению.
  Таким образом исходная задача может быть сформулирована как задача невыпуклой оптимизации в пространстве порядка $10^4$ переменных.
  При этом сложность вычисления как значения функции, так и градиента, квадратична по числу переменных.
  В статье приводится сопоставительный анализ безградиентных методов с линейкой методов градиентного типа (градиентный спуск, быстрый градиентный спуск, метод сопряженных градиентов, квазиньютоновские методы) в GPU-реализациях.

\textbf{Ключевые слова:} минимизация энергии, гомологичный фолдинг, быстрый градиентный спуск, метод сопряженных градиентов, LBFGS, параллельные вычисления, GPU.
\end{abstract}

\section{Введение}

Развитие экспериментальных методов измерений в молекулярных и клеточных системах в последние десятилетия привело к значительному росту понимания биологических процессов, протекающих при различных патологических явлениях.
В свою очередь это позволяет выявлять новые мишени для проведения эффективной лекарственной терапии даже для самых тяжелых заболеваний из области онкологии и аутоиммунных нарушений.
Терапевтическими агентами в этом случае должны служить искусственно созданные биологические молекулы, обеспечивающие заданные функции.
Обеспечение возможности рационального дизайна таких молекул, а также предсказание их свойств без проведения реального биологического эксперимента, является одной из важнейших целей современной вычислительной биологии.

К сожалению, как недостаточные вычислительные мощности, так и фундаментальные ограничения используемых методов математического моделирования биологических систем не позволяют решать эту задачу полностью \textit{de novo} в общем случае.
Для этого активно развиваются методы предсказания структуры синтетических белков, имеющие, по большей части, нефизичные основания.
Другой не менее перспективной областью является модификация природных белков с целью улучшения их характеристик или придания альтернативных свойств.
Оба направления подразумевают оптимизацию геометрии исследуемых молекул с целью поиска конформации (взаимного положения всех атомов), наиболее точно соответствующей реальному состоянию биологического объекта в природе.
То, что исследуемые объекты существуют при температурах значительно превышающих ноль по шкале Кельвина, определяет постоянное движение всех компонентнов молекул, вызываемых внешними возмущениями среды и, как следствие, отсутствия равновесия в действующих силах.
Тем не менее практически полезно рассматривать некоторую статистически наиболее часто встречающуюся конформацию, как ту, в которой биологическая молекула проводит большую часть времени, а значит, скорее всего именно в ней способна выполнять свою основную функцию.
Эта конформация логичным образом будет достигаться в минимуме потенциальной энергии молекулы, определяющейся множеством попарных взаимодействий между атомами.
Природа этих взаимодействий, по большей части, лежит в области электростатики и выражается соответствующими законами, которые будет подробнее рассмотрены в п. \ref{PF} данной статьи.

Важной особенностью рассмотренного процесса является то, что белок не стремится к глобальному минимуму потенциальной энергии.
Действительно, если бы простейший белок из $100$ аминокислот начал бы сворачиваться из линейной конформации в соответствии с поиском глобального минимума, то ему пришлось бы перебирать примерно $10^{100}$ возможных конформаций, что привело бы ко времени сворачивания порядка $10^{80}$ лет даже при скоросте перебора в $10^{13}$ конформаций в секунду.
Это привело бы к тому, что за время существования нашей вселенной не свернулся бы ни один белок, что, очевидно, расходится с нашим эмпирическим опытом.
Подобные соображения были высказаны впервые в \cite{levinthal1969} и получили название "парадокса Левинталя".
Разгадкой представленного парадокса является то, что белок образуется не сразу весь а последовательно от первой аминокислоты к последней, а значит и сворачивается не весь целиком одномоментно.
Это же приводит нас к соображениям, что локально схожие последовательности аминокислот будут давать локально схожие структуры.
На основе этой идеи разработаны различные методы, позволяющие получать неплохие приближения пространственных конформаций белков с использованием баз уже известных структур.
Сами эти методы находятся за рамками рассмотрения данной статьи, но упомянуты для следующего важного соображения: в проблемах вычислительной биологии важной задачей становится локальная минимизация энергетического потенциала, не вызывающая значительного возмущения геометрии рассматриваемой биологической молекулы.

Другой интресующей вычислительных биологов проблемой является поиск места состыковки двух взаимодействующих белков.
Данный процесс подразумевает переход молекул в альтернативные конформации, энергетически более выгодные для взаимодействия.
Здесь вновь возникает задача минимизации энергии, и вновь она требует лишь локальной минимизации, поскольку поиск глобального минимума означал бы значительную "поломку" имеющейся структуры невзаимодействующих белков, что сразу нарушало бы физичность нашей модели.
 
Именно эти соображения стали мотивацией к данной работе, рассматривающей с точки зрения современных численных методов оптимизации задачу минимизации потенциальной энергии пространственной структуры белка в предположении, что имеется достаточно хорошие начальные приближения к нужному решению (локальному минимуму).

Структура статьи следующая. В п. \ref{PF} приводится постановка задачи, как задачи невыпуклой оптимизации. В п. \ref{Disc} в хронологическом порядке приводится эволюция понимания задачи. Взрывное поведение минимизируемого потенциала при сближении атомов подсказывает, что градиентные методы должны тут плохо работать. А структура функционала позволяет эффективно пересчитывать значение энергии при шевелении одного лишь атома. Эти наблюдения приводят к безградиентному методу со случайным выбором атома для шевеления на каждом шаге. Тем не менее, сравнительный анализ этого метода и градиентного спуска показал, что при наличии GPU на практике более эффективным оказалось использование градиентного спуска с параллелизацией вычислений на GPU (для градиентного спуска можно в полной мере использовать потенциал современных видеокарт с тысячами вычислительных ядер). В п. \ref{CG} развивается стандартная сюжетная линия в невыпуклой оптимизации, которая сейчас активно эксплуатируется в обучении глубоких нейронных сетей \cite{Goodfellow2017}. Эта линия заключается в том, что если градиентные методы неплохо сходятся, то следует пробовать сделать следующий шаг -- рассмотреть ускоренные (быстрые) градиентные методы, которые доказуемо дают ускорение только для выпуклых задач. А дальше можно рассмотреть и различные варианты методов сопряженных градиентов и квази-ньютоновских методов, гарантированно сходящиеся только для выпуклых задач. В данном пункте приводится современный обзор результатов по скорости сходимости данных методов для выпуклых задач. Обсуждаются также методы более высокого порядка, к которым можно отнести, например, метод Ньютона.
В п. \ref{Practice} приведены результаты многочисленных экспериментов с методами, описанными в предыдущем пункте. В целом, оценки скорости сходимости (и опыт практического использования) рассмотренных методов в выпуклом случае неплохо перенесся для данного класса задач и на невыпуклый случай, что изначально было совершенно неочевидно в виду наличия большого числа взрывных особенностей у оптимизируемого функционала. Наиболее эффективным оказались быстрый градиентный метод Нестерова \cite{nesterov2013introductory} с специальной стратегией выбора шага, квази-ньютоновский метод LBFGS с памятью 3 итерации и один из вариантов метода сопряженных градиентов Полака--Рибьера--Поляка \cite{nocedal2006sequential}. Все методы программировались сначала на Haskell (создавался промышленный код; использование функционального языка практически исключало наличие возможных ошибок -- при их наличии не удалось бы скомпилировать программу) с использованием библиотеки Accelerate \cite{Haskel} для работы с архитектурой GPU, а затем на CUDA С также для работы c GPU. К сожалению, использование библиотеки Accelerate заметно (в несколько раз) дополнительно замедляло время работы программы и привносило дополнительные сложности при запуске программы на имеющимся сервере по сравнению с более стандартной в данном контексте архитектурой CUDA C \cite{zhmur}. Поэтому итоговый промышленный продукт, который сейчас активно используется в различных сервисах компании БИОКАД, был написан на CUDA C. Некоторые особенности реализации отмеченных методов (например, способ выбора шага/осуществление одномерной вспомогательной минимизации) также обсуждаются в данном пункте.

В данной работе авторы и компания БИОКАД, ставшая инициатором проводимых исследований, совершенно намерено пошли на изложение всех основных деталей реализации промышленного кода, которые могут представлять интерес для коллег.
В статье также описан как положительный опыт, так и отрицательный, что не является типичным для научных статей.
\section{Постановка задачи} \label{PF}

Минимизируемая геометрия белка будет представляться координатами его атомов.
Для каждой из двадцати видов аминокислот, используемых при построении белков в живой природе, можно определить конкретное количество атомов, в нее входящих.
В среднем это число будет составлять около двух десятков.
Каждый атом, в свою очередь, имеет три пространтсвенные коориданты $\{x_k ,y_k ,z_k\}$.
Далее мы будем работать в этом координатном пространстве.
Обозначим через $n_{ }\left(n\sim 10^4\right)$ размерность этого пространства.
Отметим, что в литературе часто работают и в другом пространстве: в пространстве углов и расстояний между соседними атомами, см., например, \cite{canutescu2003cyclic}, \cite{coutsias2004kinematic}. 

Центральное место в задаче минимизации потенциальной энергии молекулярного комплекса, разумеется, занимает приближение межатомных взаимодействий с помощью классических потенциальных сил.
Различные функциональные формы, а также используемые в данных формулах коэффициенты принято называть \textit{силовыми полями} \cite{Ponder2003}.
Получение таких полей -- тяжелая и кропотливая работа, включающая в себя многократную постановку сложных физических экспериментов и проведение квантово-механических вычислений.
В связи с этим коэффициенты силовых полей могут разниться и с разным качеством описывать различные виды молекулярных систем.
Тем не менее для белков, ввиду малого количества видов аминокислот (по сравнению с общим пространством допустимых органических соединений), применимы многие из этих полей: AMBER \cite{Cornell1995}, CHARMM36 \cite{Huang2017}, OPLS \cite{Jorgensen1996} и другие.

Поле OPLS, особенно в последних редакциях, отлично зарекомендовало себя в точных оценках потенциальной энергии \cite{Sweere2017}, в связи с чем именно оно было
выбрано в качестве основной рабочей модели.
Разберем его устройство.

В состав поля OPLS входят два вида основных взаимодействий: энергия связей и энергия несвязанных атомов.
Первый вид энергетической функции рассматривает каждую связь между парой атомов как некоторый протяженный трехмерный объект, способный растягиваться (stretch), гнуться (bend) и скручиваться (torsion).
Каждый из видов деформации приводит к изменению потенциальной энергии связей в соответствии с законами упругой деформации.

Так, энергию растяжения можно получить, просуммировав потенциальную энергию по всем парам связанных атомов:
$$
E_{stretch} = \sum\limits_{bonds}K(r - r_0)^2
$$
Текущая длина связи представлена через $r$ и может быть легко вычислена путем вычитания координат и взятия нормы вектора.
Помимо этого параметра в данной формуле присутствуют два коэффициента: $K$ (коэффициент растяжения связи) и $r_0$ (равновесная длина связи).
Коэффициенты являются табличными константами, определенными для каждой пары атомов, в соответствии с их типом.
Например: $(N, C_{\alpha}), (C, N), (C_{\alpha}, C_{\beta}^{(Ala)})$ и т.д.

Аналогично, используя коэффициент изгиба $K_{\theta}$ и равновесный угол $\theta_0$ можно ввести энергию изгиба, пользуясь всеми тройками последовательно связанных атомов:
$$
E_{bend} = \sum\limits_{angles}K_{\theta}(\theta - \theta_0)^2
$$

Потенциальная энергия кручения определяется через двугранный угол $\phi$ между плоскостями, заданных осью, проходящей через крутящуюся связь и двум атомам, прилегающим с концам связи с разных сторон.
Так как кручение связи имеет периодический характер, то и энергия также будет описываться некоторым периодическим законом.
В поле OPLS все четверки последовательно связанных атомов образуют последний компонент связанной энергии следующим образом:
$$
\begin{aligned}
  E_{torsion} = \sum\limits_{dihedrals}\frac{1}{2}(V_1(1 + cos(\phi)) &+ V_2(1 - cos(2\phi))\\
  &+ V_3(1 + cos(3\phi)) + V_4(1 - cos(4\phi)))
\end{aligned}
$$
Для вычисления потребуется 4 табличных параметра: $V_1$, $V_2$, $V_3$, $V_4$.

Три рассмотренных вида энергии дают нам полное представление обо всем потенциале связанных атомов:
$$
E_{bonded} = E_{stretch} + E_{bend} + E_{torsion}
$$
Помимо этой энергии все пары атомов взаимодействуют благодаря электростатическим силам, действующих вне зависимости от наличия связи.
Такие силы раскладываются на две основных составляющих: кулоновские и силы Ван-дер-Ваальса, включающие в себя три вида слабых электромагнитных взаимодействий.
Закон Кулон в хорошо известной форме полностью ложится в силовое поле OPLS:
$$
\begin{aligned}
E_{cul} &= \sum\limits_{i < j} \frac{q_i q_j e}{4 \pi \epsilon_0 r_{ij}} \\
       &= C \cdot \sum\limits_{i < j} \frac{q_i q_j}{r_{ij}},\; \text{где } C = \frac{e}{4 \pi \epsilon_0} = 1389.38757  
\end{aligned}
$$
Здесь заряды атомов $q_i$ являются табличными константами, а расстояние, как обычно, параметром.

Ван-дер-Ваальсовы силы возникают вследствии поляризации атомов в присутствии друг друга и описываются известным потенциалом Леннарда-Джонса (потенциал $6-12$), который задает экспоненциальное отталкивание на расстояниях меньше радиуса атома и притяжение на расстояниях больше.
Коэффициентами функции являются Ван-дер-Ваальсовый радиус атома $\sigma$ и глубина потенциальной ямы $\varepsilon$, параметром опять же расстояние между атомами:
$$
E_{vdw} = 4 \cdot \sum\limits_{i < j} \sqrt{\varepsilon_i \varepsilon_j} \left( \left( \frac{\sqrt{\sigma_i \sigma_j}}{r_{ij}} \right)^{12} - \left( \frac{\sqrt{\sigma_i \sigma_j}}{r_{ij}} \right)^{6} \right)
$$

Итого, задача минимизация потенциала полной модел и силового поля OPLS будет выглядеть следующим образом:
$$
E(\{r\},\{\theta\},\{\phi\}) = E_{stretch} + E_{bend} + E_{torsion} + E_{cul} + E_{vdw} \to \min_{\{r\},\{\theta\},\{\phi\}}.
$$
Заметим, что набор переменных, по которым происходит оптимизация $\{r\}$, $\{\theta\}$, $\{\phi\}$ однозачно определяется положениями атомов $\{x_k ,y_k ,z_k\}$. Как уже отмечалось, именно в пространстве $\{x_k ,y_k ,z_k\}$ далее будет решаться выписанная задача. Отметим также, что значения коээфициентов (т.е. всех параметров, не являющихся переменными, по которым осуществляется оптимизация: $\{r\}$, $\{\theta\}$, $\{\phi\}$) определяются структурой рассматриваемой макромолекулы, и на данный момент не известны какие-то красивые (компактные) способы их задания, кроме как табличным образом \cite{Jorgensen1996}.

Для всех описанных выше функций $E_{bonds}$, $E_{bend}$, $E_{dihedrals}$ можно посчитать значение функции и градиента за ${\rm O}\left( n \right)$ арифметических операций, а для $E_{cul}$, $E_{vdw}$ -- за ${\rm O}\left( n^2 \right)$.
При этом результат о сложности вычислении градиента можно получить с помощью автоматического дифференцирования \cite{nocedal2006sequential}.
Однако, в данном случае, лучше выписать явные формулы для расчета градиента, не требующие выгрузки в (оперативную) память дерева вычислений.

Важно отметить, что рассматривая белки в прикладных задачах, мы всегда подразумеваем наличия некоторой внешней среды, в которой они присутствуют.
Обычно эта среда представляет из себя раствор электролитов, то есть молекулы воды с ионами тех или иных солей.
Представлять подобный раствор можно множеством молекул, взаимодействующих с белком по тем же законам.
Однако такое представление оказывается вычислительно сложным и излишним в задаче минимизации энергии белка, так что вместо явного представления растворителя можно обратиться к приближению непрерывной среды.
Данное приближение дает возможность пользоваться некоторым потенциалом силы $E_{solvation}$, действующей между растворителем и атомами рассматриваемых макромолекул.
Можно выписать расчетные формулы для $E_{solvation}$ \cite{still1990semianalytical}, однако они достаточно громоздки, и могут быть информативными только для специалистов, поэтому было решено их здесь не приводить.
Далее будет важно только то, что функция $E_{solvation}$ локально достаточно гладкая (можно посчитать градиент), а сложность вычисления функции и ее градиента -- ${\rm O}\left( n^2 \right)$, что соответствует (по порядку) сложности вычисления всех остальных слагаемых.
В ряде экспериментов вычисления проводились без этого слагаемого.

Задача, которую необходимо решать, может быть сформулирована следующим образом: найти в введенном координатном пространстве такую конфигурацию, которая бы минимизировала функцию $E$, если известно достаточно хорошее начальное приближение.
Именно последнее предположение отличает эту статью от подавляющего большинства других работ, посвященных белковому фолдингу, см., например, \cite{zhmur}. 

Подчеркнем, что не требуется искать глобальный минимум, как было сказано ранее, это будет неверно с физической точки зрения.
Нашей же задачей станет поиск именно локального минимума, в который можно попасть из рассматриваемой точки старта.
Отметим при этом, что используя различные обобщения конструкции работы \cite{wales1997global} (например, Monotonic Basin Hopping -- MBH \cite{posypkin2010}), можно пытаться строить на базе описываемых далее локально сходящихся методов, новые методы, которые находят более глубокий минимум.
Однако стоит отметить, что в подавляющем числе проведенных экспериментов не удалось сколько-нибудь существенно (больше, чем на 10\% от абсолютного значению функционала в найденном локальном минимуме) улучшать найденный локальный минимум.
При этом процедура MBH запускалась на GPU на полный день, в то время как локальная минимизация работала десятки секунд.
\section{Обсуждение возможных подходов к численному решению задачи}\label{Disc}

\begin{figure}
\begin{center}
\includegraphics[width=10cm]{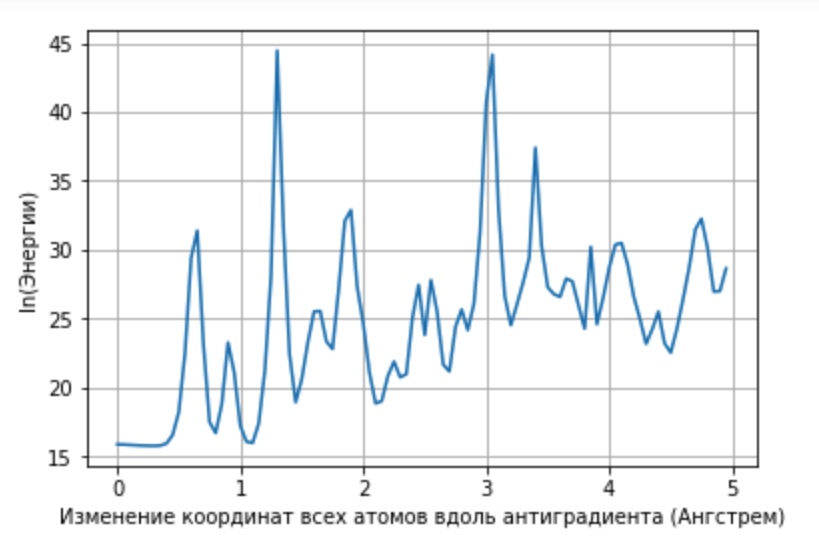}
\end{center}
 \caption{Существенная невыпуклость задачи минимизации OPLS force field}
  \label{fig1D}
 \end{figure}

Прежде всего заметим, что рассматриваемая задача является существенно невыпуклой. Особенно в окрестности начальной конфигурации. Типично, при выставлении начальной конфигурации (позы), некоторые атомы оказываются неправдоподобно близко к друг другу, что приводит к огромным значениям начальной энергии ($\sim 10^{8}$ кДж/моль вместо $\sim -10^4$ кДж/моль в минимуме). Приведем один характерный пример. Если зафиксировать единичный вектор, направленный, вдоль антиградиента функции $E$, вычисленного в начальный момент, и двигаться вдоль этого направления, то можно наблюдать следующее поведение энергии $E$, см. рис. \ref{fig1D}.

Более того, из структуры оптимизируемого функционала и рис. \ref{fig1D} понятно также, что оптимизируемый функционал может рассматриваться как гладкий только локально. Константа Липшица (и, тем более, константа Липшица градиента) не является равномерно ограниченной в достаточно большой, но ограниченной окрестности точки старта. Но в таком случае для минимизации такого функционала нельзя использовать градиентные методы -- нет гарантий даже локальной сходимости \cite{ghadimi2015generalized}, \cite{nesterov2017random}. 
 
Для решения таких задач глобальной оптимизации теория рекомендует использовать безградиентные методы типа simulated annealing или марковского поиска \cite{zhigljavsky2007stochastic}, \cite{zhigljavsky2012theory}. Ниже описывается подход, который был выбран изначально для решения задачи в отсутствии возможности использовать распараллеливание на GPU, но в условиях возможности использовать распараллеливание на CPU. Также предполагалось отсутствие составляющей  $E_{solvation}$ в $E$. В качестве языка программирования для реализации описываемого далее (безградиентного) метода использовался (интерпретируемый) язык Python 3.
 
\begin{figure}
\begin{center}
\includegraphics[width=.5\linewidth]{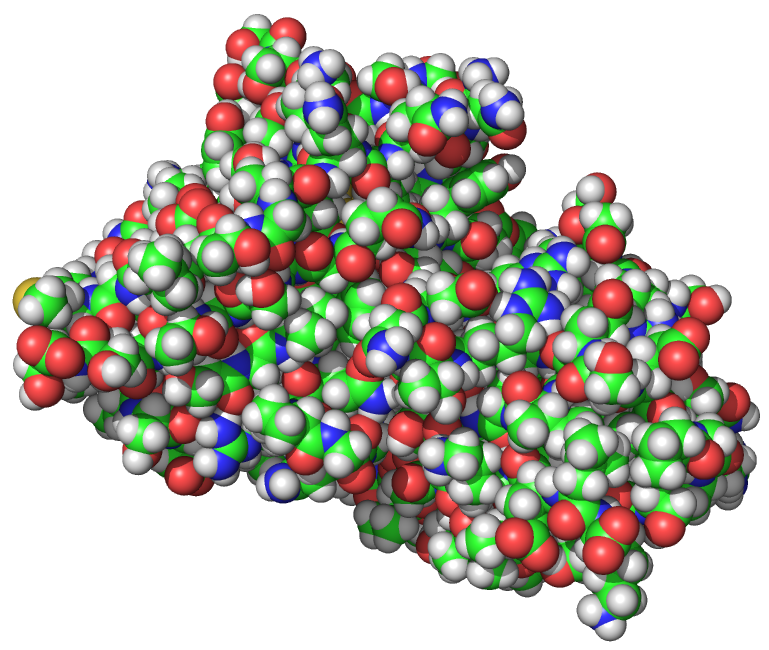}
\end{center}
\caption{Структура типичного белка}
\label{belok}
\end{figure}

В основе предлагаемого подхода -- ``шевеление'' на каждой итерации только 
одного (случайно выбранного) атома при ``замороженных'' остальных. Заметим, 
что при изменении положения одного атома пересчет $E_{stretch} $, $E_{bend} 
$, $E_{torsion} $ будет стоить ${\rm O}\left( 1 \right)$, поскольку 
затрагивает только ``соседние'' по химическим связям атомы (``соседние'' 
взяты в кавычки, потому что речь идет не только о непосредственных соседях, 
но и о соседях через две и даже три химические связи). Обычно число таких 
``соседей'' не больше 15. Для приближенного пересчета $E_{vdw} $ можно использовать 
специальные структуры данных типа Kd Tree \cite{de2000computational}, которые имеются, например, в 
библиотеке SciPy Python 3. Однако оптимально использовать в данном случае 
деревья диапазонов (Range Tree), описанные, например, в главе 5 \cite{de2000computational}. С помощью такой структуры данных первый пересчет $E_{vdw} $ будет занимать 
${\rm O}\left( {\ln ^2n} \right)$ -- это время уходит на поиск ближайших 
пространственных соседей рассматриваемого атома. Последующие пересчеты 
$E_{vdw} $ при изменении положения того же атома будут занимать ${\rm 
O}\left( 1 \right)$. К сожалению, для небольших белков описанный способ пересчета $E_{vdw}$ не дает на практике ускорения по сравнению с честным пересчетом. Связано это с тем, что в число пространственных соседей, которых надо учитывать, для рассматриваемого атома попадают обычно десятки (а иногда и сотни) других атомов, см. рис. \ref{belok} (иначе точность аппроксимации будет не достаточной). Еще сложнее дело обстоит с $E_{cul} $. 

Численные эксперименты убедительно свидетельствовали о том, что для пересчета значений $E_{cul} $ необходимо учитывать все слагаемые в сумме, которые связаны с рассматриваемым атомом. Действительно, относительный вклад в потенциал неучтенных частей $E_{vdw}$ для шара радиуса 7 A (радиус атома водорода 1 А -- один ангстрем) с центром в данном атоме составляет $\sim 10^{-5}$. Другое дело для электростатической составляющей $E_{cul}$ -- для этого же шара относительный вклад будет $\sim 10^{-3}$, что уже сопоставимо с возможной энергетической выгодой, которая получается при шевелении атома. Иначе говоря, нельзя не учитывать в электростатическом потенциале слагаемые вне шара! Поскольку всего таких слагаемых 
${\rm O}\left( n \right)$, то получается, что сложность пересчета $E_{cul} $ будет ${\rm 
O}\left( n \right)$, что на два порядка больше, чем по остальным 
компонентам энергии. Далее излагается  оригинальная конструкция, позволяющая получить амортизационную сложность пересчета 
$E_{cul} $ порядка ${\rm O}\left( 1 \right)$. 

Прежде всего, заметим, что на одной итерации необходимо посчитать (точнее, 
пересчитать) значение полной энергии при изменении положения одного атома ни 
один раз, а много раз, чтобы определить наилучшее положение этого атома при 
фиксированных положениях всех остальных атомов. Выберем число таких 
пересчетов порядка $n$. Покажем, как можно осуществить все эти $\sim n$ 
пересчетов за время ${\rm O}\left( n \right)$. Пусть на рассматриваемой 
итерации был выбран атом с номером $i$. Обозначим приращения положения этого 
атома через $\{\delta _i^x ,\delta _i^y ,\delta _i^z \}$. Тогда за время 
${\rm O}\left( {\ln ^2n} \right)$ (см. выше) сумму $E_{cul}$ можно разбить на две 
неравнозначные составляющие. Первая составляющая (${\rm O}\left( 1 \right)$ 
слагаемых) будет содержать сумму по ближайшим пространственным соседям атома 
$i$, а вторая сумма будет содержать все, что не вошло в первую (${\rm 
O}\left( n \right)$ слагаемых). Вторую сумму с хорошей точностью можно 
заменить ее рядом Тейлора по $\{\delta _i^x ,\delta _i^y ,\delta _i^z \}$ с 
членами до первого порядка включительно. 
Несложно понять, что вторая сумма есть просто линейная форма 
$C_x \delta _i^x + C_y \delta _i^y +C_z \delta _i^z $, коэффициенты которой 
$\left\{ {C_x ,C_y ,C_z } \right\}$ достаточно посчитать один раз (это можно 
сделать за ${\rm O}\left( n \right))$. Таким образом, чтобы посчитать 
значение $E_{cul}^i \left( {\delta _i^x ,\delta _i^y ,\delta _i^z } 
\right)$ (измнение энергии $E_{cul}$, возникшее вследствии шевеления атома $i$ на ${\delta _i^x ,\delta _i^y ,\delta _i^z }$) один раз нужно время ${\rm O}\left( n \right)$, но чтобы посчитать 
$E_{cul}^i \left( {\delta _i^x ,\delta _i^y ,\delta _i^z } \right)$ порядка 
$n$ раз также нужно время ${\rm O}\left( n \right)$. Получается, что среднее 
время на один расчет будет ${\rm O}\left( 1 \right)$. Таким образом, 
пересчитать полную энергию порядка $n$ раз при изменении положения одного 
атома (случайно выбранного на данной итерации) можно за время ${\rm O}\left( 
n \right)$. Это и будет стоимостью итерации, на которой приближенно ищется 
оптимальное положение случайно выбранного атома. Ясно, что число таких 
итераций (в виду их дешевизны) можно сделать очень большим. Более того, если 
в рассматриваемой эпохе случайно ``нарезать'' исходную молекулу на 
приблизительно одинаковые части и выбирать независимо в каждой части атомы 
для шевеления (вдали от границ), то можно организовать параллельный запуск 
описанного выше метода. В конце каждой эпохи все части снова собираются 
вместе, и происходит новое (случайное и независимое) нарезание на части, 
начинается новая эпоха {\ldots} Отметим, что для осуществления такого 
подхода весьма существенно, что атомы для шевеления выбираются случайно. 
Впрочем, немного более громоздкие рассуждения, возможно, позволят ускорить 
сходимость метода за счет выбора для шевеления не случайного атома, а атома, 
от шевеления которого может быть наибольший эффект -- правило 
Гаусса--Саузвелла \cite{nutini2015coordinate} (соответственно, при параллельном подходе нужно будет 
выбирать такой атом среди атомов, принадлежащих рассматриваемой части). К 
сожалению, последнее (чувствительность к шевелению) можно определить лишь 
локально -- по пересчету частных производных. Тем не менее, подобно 
описанному выше пересчету значений функции можно показать, что с аналогичной 
сложностью можно пересчитывать и частные производные, поскольку пересчитывать их 
нужно будет только у ``соседних'' по химически связям и пространственно 
атомов. По-видимому, это общий факт, связанный с развитием идей 
автоматического дифференцирования \cite{nocedal2006sequential}. Используя далее, скажем, бинарные 
кучи \cite{cormen2001introduction} можно поддерживать все время быстрый доступ к нужной для работы 
алгоритма максимальной компоненте градиента.

Описанный выше подход был реализован. К сожалению, для белков средних размеров (около 100 аминокислот) при приближении к локальному минимуму было обнаружено, что подход перестает работать -- заметно не доходя до минимума, метод перестает улучшать честно посчитанное значение энергии. Аккуратный разбор почему это происходит, показал, что неточность аппроксимации рядом Тейлора даже при разбиение суммы $ E_{cul}$ на две, более менее, одинаковые (по числу слагаемых) подсуммы имеет порядок относительной точности $10^{-4}$ -- $10^{-3}$ сопоставимый с возможным энергетическим выигрышем, получаемым при шевелении выбранного атома. Можно, конечно, использовать квадратичную аппроксимацию, но тут уже препроцессинг (вычисление коэффициентов квадратичной формы) будет довольно затратный (не окупается). 

\begin{figure}
\begin{center}
\includegraphics[width=10cm]{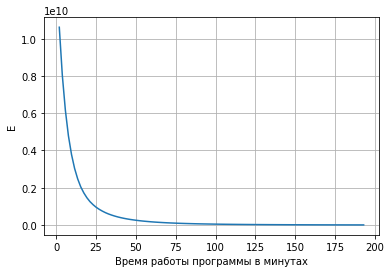}
\end{center}
\caption{Сходимость безградиентного метода (шевеление случайного атома)}
 \label{GF}
\end{figure}

Последующие эксперименты показали, что нет необходимости пересчитывать значение энергии в ${\rm O}\left( n \right)$ точках. Рассматривались следующие альтернативные стратегии: 1) выбрать случайно одну из координатных осей и посчитать значение энергии еще в двух точках на этой оси, по трем точкам построить параболу и сместить атом вдоль выбранной оси в точку минимума параболы (метод парабол); 2) аналогичный метод, но в 3D; 3) по значениям энергии в текущей точке и близкой к ней соседней точке можно оценить частную производную энергии по данной координате и далее использовать различные варианты рандомизированных покомпонентных спусков 
\cite{conn2009introduction},
\cite{ghadimi2013stochastic}, \cite{wright2015coordinate}.
Эксперименты показали, что такое огрубление (неточный поиск положения атома, который был выбран для шевеления) принципиально не менят скорости сходимости. На рис. \ref{GF} показана работа не параллельного варианта обычного покомпонентного метода с постоянным (предварительно подобранным) шагом \cite{ghadimi2013stochastic}, запущенный на современном ноутбуке (процессор 2.2 GHz; оперативная память 16 ГБ). Сделанное огрубление заметно уменьшает стоимость итерации, не сильно проигрывая в числе итераций. Итоговый общий выигрыш во времени работы обусловлен тем, что выписанные ранее оценки не учитывали числовые константы, которые в подходе с аппроксимацией $E_{cul}$ оказались достаточно большими. 

Наблюдаемая в экспериментах хорошая работа покомпонентного типа методов наводит на мысль, что и полноградиентные методы также могут хорошо работать несмотря на взрывные особенности функционала. Это, на самом деле, оказалось так. Обычный градиентный спуск в среднем за 400 итераций сходится к такому же по качеству (с точки зрения значения энергии) решению -- см. рис. \ref{GF}, затрачивая на это приблизительно такое же время (5--10 часов).

\begin{figure}
\begin{center}
\includegraphics[width=10cm]{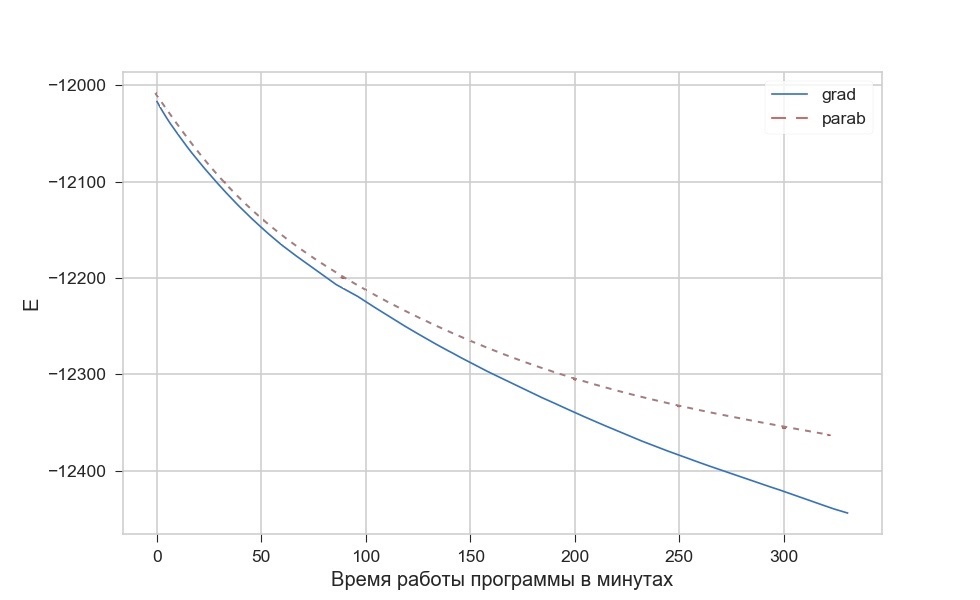}
\end{center}
\caption{Сравнение скорости сходимости безградиентного метода (шевеление случайного атома) с градиентным спуском}
 \label{Compare}
\end{figure}

На рис. \ref{Compare} и приведено сравнение времени работы безгрдиентного метода (варианта покомпонентного спуска) с градиентным спуском. Была выбрана конфигурация с изначально хорошим начальным условием (аналогичные результаты наблюдались и в общем случае). Представленный эксперимент проведен на стандартном ноутбуке(Intel Core i7-3630QM 2.40 GHz, 8GB) в стандартном режиме работы. Дополнительное отключение стандартных функций с целью повышения быстродействия и тюнинг операционной системы не проводились. Использовалась однопоточная реализация обоих методов на Python.

Метод, обозначенный на графике как 'grad', представляет собой реализацию градиентного спуска с адаптивным шагом. Адаптивные стратегии выбора шага в процедурах типа градиентного спуска обсуждаются в заключительном п. 5. 

В безградиентном методе, обозначенном как 'parab', случайным образом выбирается атом в молекуле (для шевеления) и вычисляет значение энергии при 7 различных положениях этого атома: начальное положение, и по две дополнительные позиции, симметрично расположенные относительно начального положения, вдоль каждой из трех координатных осей. По этим 7 точкам строится параболическая аппроксимация целевой функции (энергии) и вычисляется положение ее минимум. В итоге получается 8 точек, из которых выбирается та, которая доставляет минимальное значение целевой функции. Эта точка и принимается за новое положение системы атомов. 

\begin{table}[h]
	\caption{Сравнение характеристик методов}
	\label{tabular:timesandtenses}
	\centering
	\resizebox{\textwidth}{!}{
		\begin{tabular}{|c|c|c|}
			\hline
			\fontsize{12pt}{12pt}\selectfont {\bfseries Показатели} & {\bfseries ``Шевеление'' атома} & {\bfseries Адаптивный градиентный спуск} \\
			\hline
			\fontsize{12pt}{12pt}\selectfont Время 1 итерации, секунды    & 1.2075 & 122.9814 \\
            \hline
			\fontsize{12pt}{12pt}\selectfont Энергия 1 итерации, кДж/моль & 0.0225 & 2.7081 \\
			\hline
			\fontsize{12pt}{12pt}\selectfont $\sim\Delta$Энергии к 300 минуте, кДж/моль & $- 347$ & $- 413$ \\
			\hline
		\end{tabular}
	}
\end{table}

Безградиентные методы, связанные с шевелением случайно выбранного атома (и их различные модификации), могут быть распараллелены на CPU, но не более чем на нескольких десятках процессоров (для белков стандартных размеров), потому что нельзя параллельно шевелить атомы, близкие к выбранному. Более того, процедура распараллеливания совсем нетривиальная. В итоге тут так и не удалось написать быстро работающую параллельную версию программы. При этом использовать в таком подходе всю мощь имеющихся GPU для ускорения расчетов по шевелению одного атома не представляет большого смысла, поскольку сложность каждого шевеления всего ${\rm O}\left( n \right)$. Совсем другое дело распараллеливание вычислений $E_{vdw}$, $E_{cul}$ (и $E_{solvation}$, если учитывается) и их градиентов на GPU. Структура этих функционалов (вида двойной суммы) с диапазоном суммирования в каждой из сумм по порядку равным числу ядер на видеокарте (использовалась видеокарта Tesla V100 c 5120 процессорами/ядрами (1.3 GHz) и общей оперативной памятью 32 ГБ) подсказывает естественный способ распараллеливания вычислений этих двойных сумм -- по внешней сумме. А именно, внешняя сумма представляет собой сумму ${\rm O}\left( n \right)$ слагаемых, каждое слагаемое, в свою очередь представляет собой сумму ${\rm O}\left( n \right)$ слагаемых. Каждую внутреннюю сумму вычисляет свой процессор (ядро) видеокарты. 

Поскольку рассматривается задача является существенно невыпуклой задачей оптимизации, то можно лишь надеяться на локальную сходимость алгоритмов. Как известно, для достаточно гладких задач (не наш случай) обычный градиентный спуск будет сходиться к локальному экстремуму таким образом, что норма градиента после $N$ итераций в общем случае будет убывать как $\sim N^{-1/2}$. Причем в классе методов первого порядка эта оценка не улучшаема  \cite{carmon2017lower2}. Кроме того, даже если использовать методы более высокого порядка, то рассчитывать на сходимость более быструю, чем $\sim N^{-1}$ все равно не приходится \cite{carmon2017lower1}. Заметим, что обычный градиентный спуск может ``застрять'' в седловой точке \cite{nesterov2013introductory}, т.е. не в локальном минимуме. Тем не менее, недавно было показано, что ``типично'' градиентный спуск сходится именно к локальному минимуму \cite{lee2017first}, \cite{lee2016gradient}. Все эти результаты, однако, мало помогают в поиске наилучшего метода среди методов градиентного типа. Как показали численные эксперименты приведенные выше оценки локальной сходимости менее информативны (полезны), чем оценки глобальной сходимости методов в выпуклом случае. Дело в том, что (типично) в некоторой окрестности (локального) минимума можно рассчитывать на то, что поведение оптимизируемой функции с хорошей точностью соответствует поведению выпуклой функции. И хотя для рассматриваемой в  данной статье задаче априорно это совсем не очевидно, тем не менее, в экспериментах локально выпуклое поведение вполне подтверждалось.  

\begin{figure}
\begin{center}
\includegraphics[width=10cm]{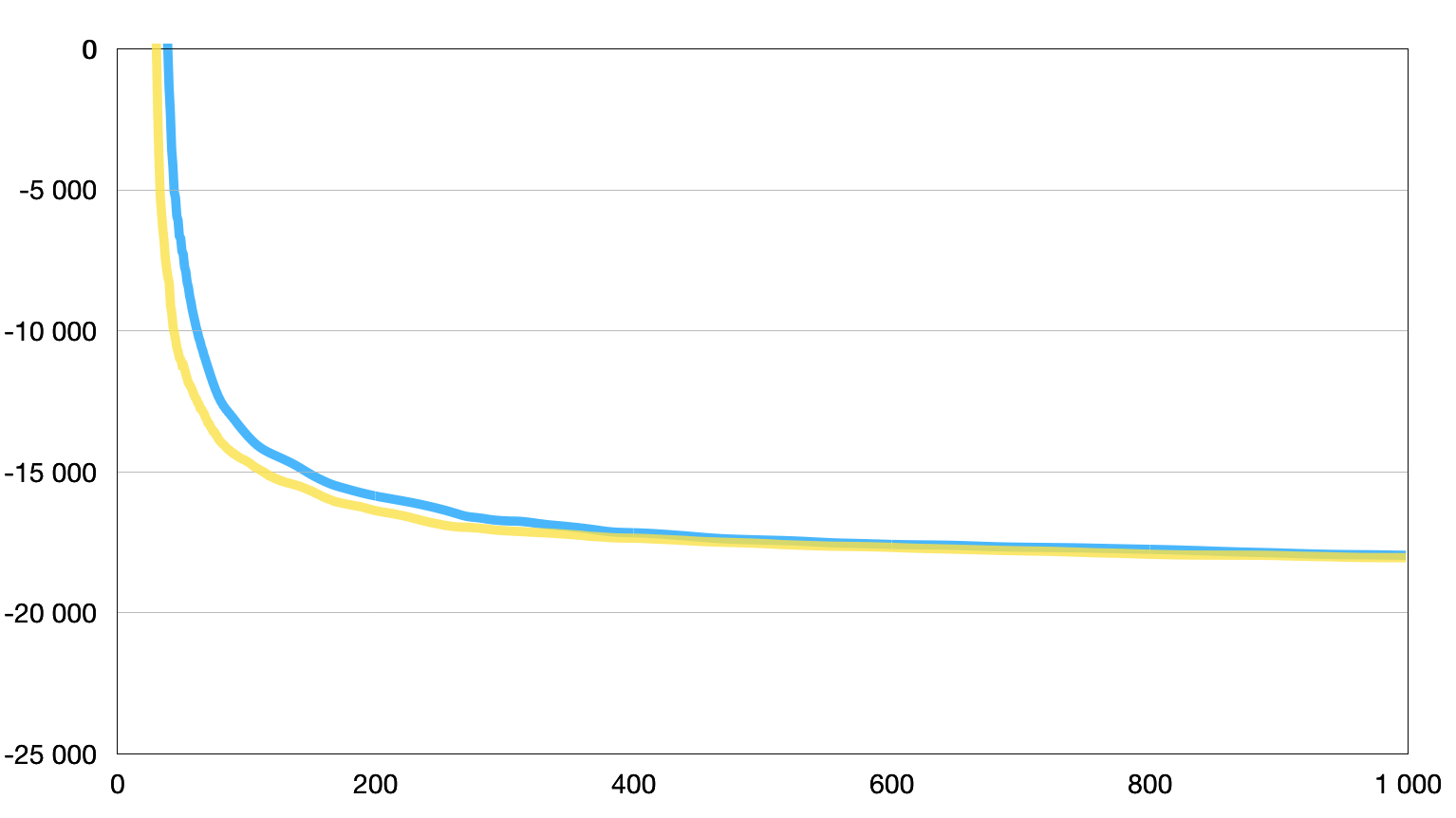}
\end{center}
\caption{Убывание энергии на итерациях градиентного спуска с правилом выбора шага по методу парабол и ускоренного градиентного спуска (см. пп. 4, 5)}
 \label{AG}
\end{figure}

Опыт использования градиентных методов для решения существенно невыпуклых задач обучения глубоких нейронных сетей \cite{Goodfellow2017} подсказывает, что если в экспериментах градиентные спуски неплохо работают, то можно попробовать использовать и ускоренные (быстрые/моментные) варианты градиентных спусков (см. п. \ref{CG}). Более того, для специальных вариантов (универсальных -- самонастраивающихся и на гладкость задачи и на ее выпуклость) ускоренных методов даже было доказано, что для невыпуклых задач они сходятся не хуже обычных градиентных методов, а в случае наличия локальной выпуклой структуры могут ускоряться, как это имеет место в глобально выпуклом случае \cite{ghadimi2015generalized}, \cite{guminov2019accelerated}, \cite{guminov2019universal}. Стоит также отметить, что ранее варианты ускоренных методов (методы типа сопряженных градиентов) уже применялись для минимизации похожих потенциалов \cite{wales1997global}. В частности, в пакете инструментов для молекулярного моделирования Schrödinger \cite{sastry2013protein} используется один из таких методов -- метод Полака--Рибьера--Поляка, см. п. \ref{CG}. Заметим также, что на рис. \ref{AG} на самом деле нарисованы два графика: график сходимости наилучшей версии (с адаптивным подбором шага по методу парабол) градиентного спуска и стандартная версия ускоренного градиентного спуска (с постоянным шагом, см. п. 4) -- ускоренный сходится немного быстрее. Естественно, в связи со всем написанным выше, попробовать пойти дальше в этом направлении и постараться подобрать наилучший ускоренный градиентный спуск, забывая на время про не выпуклость задачи. 

\section{Ускоренные градиентные методы, методы сопряженных градиентов, квазиньютоновские методы} \label{CG}

В этом пункте рассматривается общая задача безусловной выпуклой оптимизации

\begin{equation}
\label{ProbForm}
f\left( x \right) \to \mathop {\min }\limits_{x\in {\rm R}^n}.
\end{equation}

Через $x^0$ обозначим точку старта. Через $x_\ast$ -- решение задачи (\ref{ProbForm}). Если решение задачи (\ref{ProbForm}) не единственно, то под $x_*$ будем понимать такое решение задачи (\ref{ProbForm}), которое наиболее близко к $x_0$ в 2-норме. Определим $R = \left\| {x^0-x_\ast } \right\|_2 $. Будем также считать, что $\lambda_{\max} \left( \nabla^2 f\left( x \right) \right) \le L$, $\lambda_{\min} \left( \nabla^2 f\left( x \right) \right) \ge \mu \ge 0$. Введем также число обусловленности $\chi = L/\mu$. Если $n$ достаточно большое, то самым простым численным методом решения задачи (\ref{ProbForm}) будет \textit{градиентный спуск} \cite{nesterov2013introductory}, \cite{Polyak1983}:

\begin{equation}
\label{Grad}
x^{k+1} = x^k - \frac{1}{L}\nabla f\left( x^k \right).
\end{equation}
Метод (\ref{Grad}) сходится следующим образом

\begin{equation}
\label{GradConverg}
f\left( {x^N} \right)-f\left( {x_\ast } \right)\le \frac{LR^2}{2} \min \left\{ 
\frac{1}{N}, \exp\left({-\frac{N}{\chi}} \right) \right\}.
\end{equation}

В общем случае для метода (\ref{Grad}) в оценке (\ref{GradConverg}) могут быть немного улучшены лишь числовые множители (перед $L$ и $\chi$).

На практике значение параметра $L$, которое нужно методу (\ref{Grad}) для работы, как правило, не известно. К тому же, на самом деле, методу (\ref{Grad}) и не и требуется знание глобального значения этого параметра. Достаточно знать значение этого параметра на траектории метода. Возможным решением проблемы не знания $L$ является переход к \textit{методу наискорейшего спуска} \cite{Polyak1983}:

\begin{equation}
\label{GradLine}
\begin{split}
h_k \in \argmin\limits_{h\ge 0} f\left( x^k - h\nabla f\left( x^k \right) \right),\\
x^{k+1} = x^k - h_k\nabla f\left( x^k \right).
\end{split}
\end{equation}
Переход к наискорейшему спуску принципиально не меняет оценку скорости сходимости (\ref{GradConverg}).

В связи с написанным выше возникает вопрос, является ли градиентный спуск (наискорейший спуск) оптимальным методом (с точностью до числовых множителей) в классе методов вида

\begin{equation}
\label{GeneralForm}
x^{k+1}\in  x^0+\;\mbox{Lin}\left\{ 
{\nabla f\left( {x^0} \right),...,\nabla f\left( {x^k} \right)} \right\} ?
\end{equation}

Ответ на этот вопрос отрицательный \cite{Nemirovski1979}. Оказывается у градиентных (наискорейших) спусков есть различные ускоренные варианты. Поясним это на примере задачи квадратичной оптимизации.

Самой характерной задачей 
выпуклой оптимизации является задача минимизации положительно определенной 
квадратичной формы:

\[
f\left( x \right)=\frac{1}{2}\left\langle {Ax,x} \right\rangle -\left\langle 
{b,x} \right\rangle \to \mathop {\min }\limits_{x\in {\rm R}^n} .
\]

Изучив данный класс задач, можно пытаться понять, как сходятся различные методы хотя бы локально (в окрестности минимума) в задачах выпуклой оптимизации. Кроме того, такие задачи возникают как вспомогательные подзадачи при использовании методов второго порядка (например, метода Ньютона). Известно, что для задачи выпуклой квадратичной оптимизации в оценке (\ref{GeneralForm}) достаточно использовать историю только с предыдущей итерации (оптимальный метод будет находиться в этом классе). Опишем его. Это \textit{метод сопряженных градиентов }(первая итерация делается согласно (\ref{GradLine}))

\begin{equation}
\label{ConGrad}
x^{k+1} = x^k-\alpha _k \nabla f\left( {x^k} \right)+\beta _k \cdot 
\left( {x^k-x^{k-1}} \right),
\end{equation}
где

\[
\left( {\alpha _k ,\beta _k } \right) = \argmin \limits_{\alpha ,\beta } f\left( {x^k-\alpha \nabla f\left( {x^k} 
\right)+\beta \cdot \left( {x^k-x^{k-1}} \right)} \right).
\]

Метод (\ref{ConGrad}) сходится следующим (оптимальным!) образом

\begin{equation}
\label{CGConv}
\begin{split}
f\left( {x^N} \right)-f\left( {x_\ast } \right)\le \\ 
\min \left\{ 
{\frac{LR^2}{2\left( {2N+1} \right)^2},2LR^2\left( {\frac{\sqrt \chi 
-1}{\sqrt \chi +1}} \right)^{2N},\left( {\frac{\lambda _{n-N+1} -\lambda _1 
}{\lambda _{n-N+1} +\lambda _1 }} \right)^2R^2} \right\},
\end{split}
\end{equation}
где $N\le n$, $0\le \mu = \lambda_{\min}\left(A\right) = \lambda_{1} \le \lambda_{2} \le ... \le \lambda_{n} = \lambda_{\max}\left(A\right) = L$. Второй аргумент в минимуме (\ref{CGConv}) при больших значениях $N$ оценивается сверху $2LR^2\exp\left(-2N / \sqrt{\chi}\right)$. При $N=n$ метод гарантированно находит точное решение, что следует из последней оценки в минимуме (\ref{CGConv}). Сформулированный результат является 
фундаментальным фактом (жемчужиной) выпуклой оптимизации и вычислительной линейной алгебры одновременно, и базируется на наличии рекуррентных формул для \textit{многочленов Чебышёва} \cite{nocedal2006sequential}, \cite{Nemirovski1979}. 

Заметим, что можно выписать явные формулы для коэффициентов $\alpha_k$, $\beta_k$. Однако на практике более распространены варианты методов сопряженных градиентов не со вспомогательной двумерной оптимизацией или с явно выписанными шагами, а со вспомогательной одномерной оптимизацией. Наиболее популярными вариантами метода сопряженных градиентов являются следующие два метода \cite{nocedal2006sequential} (на наш взгляд  исторически правильнее называть метод Полака--Рибьера методом Полака--Рибьера--Поляка):

\[
h_k \in \argmin \limits_{h\in {\rm R}} f\left( {x^k+hp^k} \right),
\]
\[
x^{k+1}=x^k+h_k p^k,
\]
\[
p^{k+1}=-\nabla f\left( {x^{k+1}} \right) + \beta_k p^k,
\]
\[
p^0=-\nabla f\left( {x^0} \right),
\]
\[
\beta _k =-\frac{\left\| {\nabla f\left( {x^{k+1}} \right)} \right\|_2^2 
}{\left\| {\nabla f\left( {x^k} \right)} \right\|_2^2 },\quad \mbox{(формула 
Флетчера--Ривса)}
\]
\[
\beta _k =-\frac{\left\langle {\nabla f\left( {x^{k+1}} \right),\nabla 
f\left( {x^{k+1}} \right)-\nabla f\left( {x^k} \right)} \right\rangle 
}{\left\| {\nabla f\left( {x^k} \right)} \right\|_2^2 }. \quad \mbox{(формула 
Полака--Рибьера)}
\]

Для задач квадратичной оптимизации оба метода эквивалентны между собой и эквивалентны методу (\ref{ConGrad}). Для общих задач выпуклой оптимизации по этим методам не удалось пока получить оптимальные порядки скорости сходимости (установлен только сам факт глобальной сходимости для задач гладкой выпуклой оптимизации). Тем не менее, именно эти два варианта метода сопряженных градиентов наиболее часто используются при решении практических задач \cite{nocedal2006sequential} (в том числе не обязательно выпуклых). При этом с некоторой периодичностью (обычно период выбирают пропорционально размерности пространства, в котором происходит оптимизация) требуется перезапускать метод, обнуляя историю: вместо $p^{k+1}=-\nabla f\left( {x^{k+1}} \right)+\beta _k p^k$ в момент рестарта полагают $p^{k+1}=-\nabla f\left( {x^{k+1}} \right)$. По-видимому, необходимость в таких рестартах обусловлена желанием правильно сходиться в случае оптимизации сильно выпуклых функций \cite{o2015adaptive}.

Метод (\ref{ConGrad}) ничего не требует на вход (никаких параметров), а работает оптимально на классе гладких выпуклых задачи и при этом также оптимально на его подклассе -- гладких сильно выпуклых задач. Конечно, хотелось бы, чтобы и для общих задач выпуклой оптимизации метод (\ref{ConGrad}) обладал аналогичными свойствами. Однако перенести без изменений (\ref{ConGrad}) на весь класс задач выпуклой оптимизации не получилось. Тем не менее, в конце 70-х годов XX века А.С. Немировскому \cite{Nemirovski1979}, \cite{narkiss2005sequential}  удалось предложить две отдельные модификации метода (\ref{ConGrad}) для класса гладких выпуклых задач и для класса 
гладких сильно выпуклых задач, которые доказуемо сходятся (с точностью до числового множителя) по оценкам, соответствующим первому аргументу минимума в (\ref{CGConv}) и второму (в сильно выпуклом случае). Первый метод (для выпуклых задач) также не требуют на вход никаких параметров. Второй метод (для сильно выпуклых задач) использует процедуру рестартов, и в общем случае требует знания параметра сильной выпуклости. При этом оба метода также как и метод сопряженных градиентов требуют на каждой итерации решения вспомогательной малоразмерной задачи выпуклой оптимизации. Ближе к концу данного пункта будет приведен пример одного из методов такого типа (\ref{optim}).

Вспомогательные маломерные оптимизации на каждой итерации требуют дополнительных затрат. Вспоминая, что для квадратичных задач можно выписать явные выражения для параметров $\alpha_k$, $\beta_k$ (см., например, \cite{Polyak1983}), естественно, задаться вопросом: можно ли получить аналогичные по порядку оценки скорости сходимости без вспомогательной маломерной оптимизации на каждой итерации. Первым ускоренным градиентным методом с постоянными шагами (это еще более сужает класс методов) для не квадратичных задач выпуклой оптимизации был (двухшаговый) \textit{метод тяжелого шарика}

\begin{equation}
\label{HB}
x^{k+1}=x^k-\alpha \nabla f\left( {x^k} \right)+\beta \cdot \left( 
{x^k-x^{k-1}} \right),
\end{equation}
предложенный Б.Т. Поляком в 1963--1964 гг. \cite{Polyak1964}. Локальный анализ скорости его сходимости (с помощью \textit{первого метода Ляпунова}) при специальном выборе параметров шага $\alpha ,\beta >0$ давал правильные порядки локальной скорости сходимости в сильно выпуклом случае --  отвечающие с точностью до числовых множителей второму аргументу минимума в (\ref{CGConv}) \cite{Polyak1983}. Однако с установлением глобальной сходимости были некоторые трудности. В частности, на специально подобранном примере с разрывным гессианом метод может и не сходиться \cite{lessard2016analysis}, а в Чезаровском смысле траектории метода сходятся медленнее -- аналогично (неускоренному) градиентному методу \cite{ghadimi2015global}. Несмотря на отмеченные сложности метод тяжелого шарика, по-прежнему, активно используется и продолжает развиваться \cite{loizou2017linearly}. В 1982--1983 гг. Ю.Е. Нестеров в кандидатской диссертации (научным руководителем был Б.Т. Поляк) предложил первый ускоренный (быстрый) градиентный метод с постоянными шагами, для которого удалось доказать глобальную сходимость с оптимальной (с точностью до числового множителя) скоростью \cite{Nesterov1983}. Метод был ``забыт'' почти на 20 лет. Большое внимание этот метод привлек к себе лишь после выхода в 2004 году монографии \cite{nesterov2013introductory}. Важную роль в привлечении внимания к методу сыграла также работа \cite{beck2009fast}, имеющая большое число цитирований.

Для задач безусловной минимизации наиболее популярным сейчас является следующий (двухшаговый) вариант быстрого (ускоренного) градиентного метода (Нестерова) 

\[
x^{k+1}=y^k-\frac{1}{L}\nabla f\left( {y^k} \right),
\]
\[
y^k=x^k+\frac{k-1}{k+2}\left( {x^k-x^{k-1}} \right),
\]
который также можно понимать как \textit{моментный метод} \cite{su2014differential}:

\begin{equation}
\label{moment}
x^{k+1}=x^k-\frac{1}{L}\nabla f\left( {x^k+\frac{k-1}{k+2}\left( 
{x^k-x^{k-1}} \right)} \right)+\frac{k-1}{k+2}\left( {x^k-x^{k-1}} \right).
\end{equation}

Приведем также вариант моментного метода для $\mu $-сильно выпуклой функции в 2-норме \cite{nesterov2013introductory}:

\begin{equation}
\label{momentS}
x^{k+1}=x^k-\frac{1}{L}\nabla f\left( {x^k+\frac{\sqrt L -\sqrt \mu }{\sqrt 
L +\sqrt \mu }\left( {x^k-x^{k-1}} \right)} \right)+\frac{\sqrt L -\sqrt \mu 
}{\sqrt L +\sqrt \mu }\left( {x^k-x^{k-1}} \right).
\end{equation}

Оба метода (\ref{moment}), (\ref{momentS}) с точностью до числовых множителей сходятся согласно первому и второму аргументу в минимуме оценки (\ref{CGConv}). Более того, можно даже объединить оба метода (\ref{moment}), (\ref{momentS}) в один, сходящийся по оценке, наилучшей из двух отмеченных оценок \cite{nesterov2013introductory}. Однако все равно требуется априорно знать параметр $\mu$. На данный момент неизвестны такие варианты быстрых (ускоренных) градиентных методов для сильно выпуклых задач, которые бы в общем случае обходились без этого знания. Отметим при этом, что незнание параметра $L$ может быть устранено вспомогательной маломерной минимизацией. 

Из сопоставления (\ref{Grad}) и (\ref{moment}), (\ref{momentS}) легко заметить, что сложность (трудоемкость) итераций у градиентного спуска и у быстрого градиентного метода практически одинаковы, в то время как скорости сходимости отличаются очень существенно. Именно это обстоятельство и обусловило огромную популярность быстрых градиентных методов. 

Среди последних достижений в развитии ускоренных (быстрых) градиентных методов отметим, \textit{оптимизированный} вариант быстрого градиентного метода для задач безусловной оптимизации \cite{drori2018efficient}

\begin{equation}
\label{optim}
\begin{array}{c}
 y^{k+1}=\left( {1-\frac{1}{t_{k+1} }} \right)x^k+\frac{1}{t_{k+1} }x^0, \\ 
 d^{k+1}=\left( {1-\frac{1}{t_{k+1} }} \right)\nabla f\left( {x^k} 
\right)+\frac{2}{t_{k+1} }\sum\limits_{j=0}^k {t_j \nabla f\left( {x^j} 
\right)} , \\ 
 x^{k+1}=y^{k+1}-\frac{1}{L}d^{k+1}, \\ 
 \end{array}
\end{equation}
где

\[
t_{k+1} =\frac{1+\sqrt {4t_k^2 +1} }{2},
\quad
k=0,...,N-2,
\]
\[
t_N =\theta _N ,
\quad
\theta _{k+1} =\frac{1+\sqrt {8\theta _k^2 +1} }{2},
\quad
k=0,...,N-1.
\]

На классе гладких выпуклых задач метод (\ref{optim}) сходится согласно оценке 

\begin{equation}
\label{optimE}
f\left( {x^N} \right)-f\left( {x_\ast } \right)\le \frac{2LR^2}{\theta _N^2 
}\left( {\le \frac{LR^2}{N^2}} \right),
\end{equation}
которая достигается, например, на функции

\begin{equation}
\label{Example}
f\left( x \right)=\left\{ {\begin{array}{l}
 \frac{LR}{\theta _N^2 }\left\| x \right\|_2 -\frac{LR^2}{\theta _N^4 
},\;\left\| x \right\|_2 \ge \frac{R}{\theta _N^2 }, \\ 
 \frac{L}{2}\left\| x \right\|_2^2 ,\;\left\| x \right\|_2 <\frac{R}{\theta 
_N^2 }. \\ 
 \end{array}} \right.
\end{equation}
Оценка (\ref{optimE}) является точной оценкой скорости сходимости оптимальных методов вида (\ref{GeneralForm}) на классе гладких выпуклых задач. Другими словами, не существует такого метода вида (\ref{GeneralForm}) с фиксированными шагами (т.е. без вспомогательной маломерной оптимизации), который бы на всем классе задач гладкой выпуклой оптимизации сходился по оценке, лучшей чем (\ref{optimE}). Подчеркнем, что нельзя улучшить даже числовой множитель. Из метода (\ref{optim}) можно сделать метод, не требующий знания параметра $L$. Для этого в процедуре (\ref{optim}) следует заменить

\[
x^{k+1}=y^{k+1}-\frac{1}{L}d^{k+1}
\]
на (см. также (\ref{GradLine}))

\[
x^{k+1}=y^{k+1}-h_{k+1} d^{k+1},
\]
где

\[
h_{k+1} \in \argmin \limits_{h\in {\rm R}} f\left( 
{y^{k+1}-hd^{k+1}} \right).
\]
Такой метод также будет сходиться согласно оценке (\ref{optimE}). 

В общем случае на каждой итерации методов типа сопряженных градиентов (в частности, (\ref{optim})) вспомогательную задачу можно решить только приближенно. В \cite{Nemirovski1979} было установлено, что достаточно решать вспомогательную задачу с относительной точностью (по функции) $\delta ={\rm O}\left( {\varepsilon \mathord{\left/ {\vphantom {\varepsilon {N\left( \varepsilon \right)}}} \right. \kern-\nulldelimiterspace} {N\left( \varepsilon \right)}} \right)$, где $\varepsilon $ -- желаемая относительная точность (по функции) решения исходной задачи, а $N\left( \varepsilon \right)$ -- число итераций, которые делает (внешний) метод. Следовательно, вспомогательная одномерная задача может быть решена за \cite{nocedal2006sequential}

\[
{\rm O}\left( {\ln 
\left( {\delta ^{-1}} \right)} \right) =
{\rm O}\left( {\ln \left( {{N\left( \varepsilon \right)} \mathord{\left/ 
{\vphantom {{N\left( \varepsilon \right)} \varepsilon }} \right. 
\kern-\nulldelimiterspace} \varepsilon } \right)} \right)={\rm O}\left( {\ln 
\left( {\varepsilon ^{-1}} \right)} \right)
\]
обращений к оракулу (подпрограмме) за значением оптимизируемой функции.  С учетом того, что по теоретическим оценкам использование таких методов не позволяет в общем случае улучшать даже числовые множители (см. текст выше относительно оптимального метода (\ref{optim})), то, кажется, что стоит оставить методы со вспомогательной одномерной оптимизацией на каждой итерации в стороне, и ограничиться рассмотрением методов с фиксированными шагами и конечной памятью. Однако на практике типично, что различные варианты метода сопряженных градиентов (а их насчитывается уже, как минимум, несколько десятков \cite{floudas2008encyclopedia}) работают существенно быстрее ускоренных градиентных методов с постоянными шагами \cite{nesterov2018primal}. Причина такого различия связана с наличием последнего аргумента минимума в оценке (\ref{CGConv}) для методов типа сопряженных градиентов и следующим наблюдением: методы типа сопряженных градиентов сходятся также как ускоренные методы только на вырожденных (мало интересных / слишком простых для практики) примерах, см. (\ref{Example}). 

Хорошо известна (см., например, \cite{Polyak1983}) геометрическая интерпретация градиентного спуска, в основу которой положена замена исходной функции параболоидом вращения, касающимся её графика в текущей точке. Точка, доставляющая минимум параболоиду принимается за новое положение метода. Идея метода наискорейшего спуска, заключается в подборе кривизны параболоида с помощью решения вспомогательной задачи одномерной минимизации. В \textit{методе Ньютона} вместо параболоида вращения строится квадратичная аппроксимация оптимизируемой функции (на основе доступного гессиана функции, вычисленного в текущей точке), однако это приводит к необходимости решения на каждом шаге более сложной задачи -- минимизации квадратичной формы (обращения матрицы/решения системы линейных уравнений). 

Естественно, возникает идея построения какого-то ``промежуточного'' метода, с одной стороны не требующего вычисления (и, тем более, обращения) гессиана, а с другой стороны все-таки пытающегося как-то аппроксимировать гессиан, исходя из накопленной информации первого порядка (градиентов). В основу \textit{квазиньютоновских методов} положен следующий общий принцип построения квадратичной аппроксимации: квадратичная аппроксимация должна касаться графика оптимизируемой функции в текущей точке и иметь с ней одинаковые градиенты в точке с предыдущего шага (secant equation). Существует много различных способов удовлетворить этим условиям \cite{nocedal2006sequential}. У этих способов есть различные интерпретации, среди которых отметим понимание квазиньютоновских методов, как \textit{методов переменной метрики} (варианта метод сопряженных градиентов \cite{Polyak1983}). Наиболее интересными в практическом плане являются способы, которые требуют не более чем квадратичной (по размерности пространства) трудоемкости и памяти. Среди таких способов наиболее удачно себя зарекомендовал способ, приводящий в итоге к методу \textit{BFGS} \cite{nocedal2006sequential}:

\[
h_k = \argmin \limits_{h\in {\rm R}} f\left( {x^k-hH_k \nabla 
f\left( {x^k} \right)} \right),
\]
\[
x^{k+1}=x^k-h_k H_k \nabla f\left( {x^k} \right),
\]
\[
H_{k+1} =H_k +\frac{H_k \gamma _k \delta _k^T +\delta _k \gamma _k^T H_k 
}{\left\langle {H_k \gamma _k ,\gamma _k } \right\rangle }-\beta _k 
\frac{H_k \gamma _k \gamma _k^T H_k }{\left\langle {H_k \gamma _k ,\gamma _k 
} \right\rangle },
\]
где

\[
\beta _k =1+\frac{\left\langle {\gamma _k ,\delta _k } \right\rangle 
}{\left\langle {H_k \gamma _k ,\gamma _k } \right\rangle },
\quad
\gamma _k =\nabla f\left( {x^{k+1}} \right)-\nabla f\left( {x^k} \right),
\quad
\delta _k =x^{k+1}-x^k,
\quad
H_0 =I.
\]

В отличие от сопряженных градиентов (\ref{ConGrad}) в BFGS не обязательно точно осуществлять вспомогательную одномерную оптимизацию. В целом, BFGS оказался наиболее устойчивым (к вычислительным погрешностям) вариантом квазиньютоновских методов. Геометрия квазиньютоновских методов близка (но не идентична!) геометрии субградиентных методов с процедурой растяжения пространства (метод эллипсоидов, методы Шора \cite{Polyak1983}) и связана с типичным ``пилообразным'' поведением градиентных спусков в окрестности минимума для плохо обусловленных задач. Направления $\gamma _k $ и $\delta _k $ ``подсказывают'' направления растяжение/сжатие, и позволяют адаптивно улучшать обусловленность задачи, правильно аккумулируя собранную информацию в $H_k $. 

В теоретическом плане по квазиньютоновским методам на данный момент известно не так уж и много \cite{nocedal2006sequential}, \cite{nesterov2013introductory} : глобальная скорость сходимости для гладких задач выпуклой оптимизации в общем случае не выше, чем у обычных (неускоренных) градиентных методов (во всяком случае, только это пока удалось установить), а локальная скорость сходимости в случае невырожденного минимума -- \textit{сверхлинейная}, т.е. быстрее, чем линейная. 

Основным ограничением по использованию квазиньютоновских методов является необходимость в хранении и обновлении плотной квадратной матрицы $H_k $, что требует (в отличие от того, что имеет место для методов типа сопряженных градиентов) квадратичной памяти и квадратичного времени независимо от разреженности задачи. Это обстоятельство существенно ограничивает возможности по использованию таких методов для задач оптимизации с десятками тысяч переменных и более. Однако на практике используют в основном варианты таких методов \textit{с ограниченной памятью}, см., например, метод \textit{LBFGS} \cite{nocedal2006sequential}. В этом случае в памяти хранится не матрица $H_k $, а вектора, ее порождающие. Проблема, однако, тут в том, что с ростом $k$ размер этой памяти линейно растет. Поэтому обычно последовательности векторов $\left\{ {\gamma _t } \right\}$ и $\left\{ {\delta _t } \right\}$ хранят только с $m$ последних итераций ($m$ -- глубина памяти), и при этом полагают $H_{k-m} =I$. На практике $m$ часто выбирают совсем небольшим $m\simeq 3-5$. 

В заключении этого пункта приведем сопоставительный анализ методов первого порядка (градиентных методов) и методов более высокого порядка, которые могут использоваться для решения задач выпуклой оптимизации умеренных размеров ($n\le 10^4)$, в условиях отсутствия шума на классе достаточно гладких (условие $x,y\in {\rm R}^n$ можно ослабить, см., например, \cite{gasnikov2018global})
\[
\left\| {\nabla ^rf\left( y \right)-\nabla ^rf\left( x \right)} \right\|_2 
\le M_r \left\| {y-x} \right\|_2 ,
\quad
x,y\in {\rm R}^n,
\quad
M_r \le \infty ,
\quad
r=0,1,2,...,
\]
выпуклых задач. Заметим, что $\nabla ^rf\left( y \right)$ -- тензор ранга $r$. Поэтому следует пояснить, что понимается под 2-нормой в левой части данного неравенства. Ограничимся случаем $r=2$, тогда

\[
\nabla ^2f\left( x \right)=\left\| {{\partial ^2f\left( x \right)} 
\mathord{\left/ {\vphantom {{\partial ^2f\left( x \right)} {\partial x_i 
\partial x_j }}} \right. \kern-\nulldelimiterspace} {\partial x_i \partial 
x_j }} \right\|_{i,j=1}^n ,
\]
\[
 \left\| {\nabla ^2f\left( y \right)-\nabla ^2f\left( x \right)} \right\|_2 
=\mathop {\sup }\limits_{\left\| {x_1 } \right\|_2 \le 1} \;\mathop {\sup 
}\limits_{\left\| {x_2 } \right\|_2 \le 1} \left\langle {\left( {\nabla 
^2f\left( y \right)-\nabla ^2f\left( x \right)} \right)\left[ {x_1 } 
\right],x_2 } \right\rangle = 
\]
\[
 \mathop {\sup }\limits_{\left\| {x_1 } \right\|_2 \le 1} \;\mathop {\sup 
}\limits_{\left\| {x_2 } \right\|_2 \le 1} \left\langle {\left( {\nabla 
^2f\left( y \right)-\nabla ^2f\left( x \right)} \right)x_1 ,x_2 } 
\right\rangle
\]
В общем случае см. \cite{baes2009estimate}, \cite{nesterov2018implementable}. Отметим также, что при $r=0$: $\nabla ^0 f\left( x \right)=f\left( x \right)$. 

Для класса методов, у которых на каждой итерации разрешается не более чем ${\rm O}\left( 1 \right)$ раз обращаться к оракулу (подпрограмме) за $\nabla ^rf\left( x \right)$, $r\le 1$, оценка числа итераций, необходимых для достижения точности $\varepsilon $ (по функции), будет иметь вид

$${\rm O}\left(\min\left\{n\ln\left(\frac{\Delta f}{\varepsilon}\right),\frac{M_0^2 R^2}{\varepsilon^2}, \left(\frac{M_1 R^2}{\varepsilon}\right)^\frac{1}{2}\right\}\right)$$
где, как и раньше, $R=\left\| {x^0-x_\ast } \right\|_2 $, $\Delta f=f\left( {x^0} \right)-f\left( {x_\ast } \right)$. Данная оценка в общем случае не может быть улучшена даже если дополнительно известно, что, $M_2 <\infty $, $M_3 <\infty $, {\ldots} \cite{Nemirovski1979}. При этом данная оценка достигается \cite{nesterov2013introductory}, \cite{Nemirovski1979}.  Можно заметить, что описываемые выше методы (быстрые градиентные методы / методы сопряженных градиентов) как раз работают по последнему аргументу минимума в этой оценке.
Заметим, что если вместо $r=1$ имеет место $r=0$, то в приведенной оценке все аргументы минимума следует домножить на размерность пространства $n$ \cite{bayandina2017gradient}, \cite{dvurechensky2018accelerated}, \cite{Protasov1996}.

Отметим также, что у известных сейчас методов, отвечающих (с точностью до логарифмического множителя) первому аргументу минимума, достаточно дорогой является составляющая итерации, не связанная с вычислением градиента: $\gg n^2$ \cite{lee2015faster}.

Для класса методов, у которых на каждой итерации разрешается не более чем ${\rm O}\left( 1 \right)$ раз обращаться к оракулу (подпрограмме) за значениями $\nabla ^rf\left( x \right)$, $r\le p$, $p\ge 2$, оценка числа итераций, необходимых для достижения точности $\varepsilon $ (по функции), будет иметь вид \cite{gasnikov2018global}, \cite{monteiro2013accelerated}

$${\rm O}\left(\min\left\{n\ln\left(\frac{\Delta f}{\varepsilon}\right),\frac{M_0^2 R^2}{\varepsilon^2}, \left(\frac{M_1 R^2}{\varepsilon}\right)^\frac{1}{2},\ldots, \left(\frac{M_pR^{p+1}}{\varepsilon}\right)^\frac{2}{3p+1}\right\}\right).$$
Данная оценка в общем случае не может быть улучшена даже если дополнительно известно, что, $M_{p+1} <\infty $, $M_{p+2} <\infty $, {\ldots} \cite{nesterov2018implementable}, \cite{arjevani2017oracle}. Полезно отметить, что здесь, также как и для градиентных методов можно строить ``универсальные худшие в мире функции'' в классе выпуклых полиномов (от абсолютных значений линейных комбинаций переменных) степени $p+1$. 

Как видно из приведенной оценки, методы более высокого порядка ускоряют сходимость, однако при этом возникают сложные вспомогательные задачи, которые необходимо решать на каждой итерации 
\cite{gasnikov2018global}, \cite{nesterov2018implementable}, \cite{monteiro2013accelerated}. Поскольку в рассматриваемом в этой статье случае в виду невыпуклости и не равномерной гладкости нельзя рассчитывать на дополнительное ускорение, связанное с использованием методов высокого порядка (не выполняются с большим запасом необходимые для этого условия), а при этом сложность итерации гарантированно существенно увеличивается, то было решено ограничиться только методами первого порядка, к которым были отнесены, в частности, метод Полака--Рибьера--Поляка и метод LBFGS.

\section{Программная реализация, результаты вычислительных экспериментов} \label{Practice}

\subsection{Реализация методов оптимизации}

Программная реализация представленных в работе градиентных методов выполнена на языке C++11 с использованием технологии параллельного программирования Nvidia CUDA \cite{wilt_cuda_2013}. Возможность работы в режимах float / double (одинарная и двойная точности) реализована с применением механизма шаблонов C++ (C++ templates), что позволило получить унифицированный программный код, максимально эффективно использующий вычислительные возможности как центрального так и графического процессоров. Реализованный программный код скомпонован в виде разделяемой библиотеки (shared library) с соответствующим C-API (Application Programming Interface, программный интерфейс приложения), обеспечивающим возможность интеграции с другими программными модулями, применяемыми компанией БИОКАД.

Стоит отметить, что во всех алгоритмах выполняется нормирование текущего направления спуска. Эта процедура, не имеющая особого смысла с математической точки зрения, повсеместно применяется на практике т.к. существенно облегчает работу внутренних алгоритмов одномерного поиска (см. пункт \ref{LS_METHODS}), которые в изложенных ниже алгоритмах представлены как
\[
   h_k = \argmin\limits_{h} f( x^k + h r^k ) .
\]
Подобную запись следует интерпретировать именно как процедуру неточного одномерного поиска, а не ``математически честную'' минимизацию.  Тем не менее, каждая такая процедура либо находит релаксирующую точку, т.е. точку, в которой $f(x^k + h_k r^k) < f(x^k)$, либо выдают ошибку одномерного поиска, сигнализирующую, что наш градиентный метод ``зашёл в тупик'' и процедуру минимизации следует либо останавливать, либо использовать другое направление спуска.

\subsubsection{Метод FGM}
\label{FGM_DESC}

Метод FGM (Fast Gradient Method) является реализацией быстрого (ускоренного) градиентного метода Ю.Е. Нестерова (см. пункт \ref{CG}). Основным отличием от оригинальной версии является применение процедуры одномерного поиска на каждой итерации метода, т.к. значение параметра $L$ (константа Липшица) в нашем случае априори неизвестно. Особенностью метода является его немонотонность, т.е. в процессе спуска могут быть итерации, когда  $f(x^{k+1}) > f(x^k)$, поэтому необходимо дополнительно запоминать ``наилучшую точку'' $x^{\tilde{k}}$ которая и будет результатом работы метода:
\[
  K = [0, N],\quad \tilde{k} \in K,\quad \forall k \in K: f(x^{\tilde{k}}) \le f(x^k) ,
\]
где $N$ --- число выполненных итераций.

Описание реализованного метода FGM представлено в алгоритме~\ref{FGM_ALGORITHM}.

\begin{algorithm}[!h]
\caption{Метод FGM}
\label{FGM_ALGORITHM}
\begin{algorithmic}
\REQUIRE $x^0 \in \mathbb{R}^n$
\STATE $x_* \leftarrow x^0$
\STATE $x^{-1} \leftarrow x^0$
\STATE $\theta^{-1} \leftarrow 1$
\STATE $k \leftarrow 0$
\REPEAT
  \STATE $\theta_k \leftarrow 0.5\ \theta_{k-1} \left( \sqrt{ \theta_{k-1}^2 + 4} - \theta_{k-1} \right)$
  \STATE $\beta_k \leftarrow \theta_{k-1} (1 - \theta_{k-1}) / (\theta_{k-1}^2 + \theta_k)$
  \STATE $w^k \leftarrow x^k + \beta_k (x^{k} - x^{k-1})$
  \STATE $r^k \leftarrow -\nabla f(w^k) / \| \nabla f(w^k) \|_2$
  \STATE $h_k \leftarrow \argmin\limits_{h \ge 0} f( w^k + h r^k )$
  \STATE $x^{k+1} \leftarrow w^k + h_k r^k$
  \IF{$f(x^{k+1}) < f(x_*)$}
    \STATE $x_* \leftarrow x^{k+1}$
  \ENDIF
  \STATE Удаляем вектор $x^{k-1}$ из памяти
  \STATE $k = k + 1$
\UNTIL \textbf{convergence}.
\RETURN $x_*$
\end{algorithmic}
\end{algorithm}

\subsubsection{Метод CG}
\label{CG_DESC}

Метод CG (Conjugate Gradient method) является реализацией метода сопряжённых градиентов (см. пункт \ref{CG}). По умолчанию используется вариант Полак--Рибьер--Поляк (как показавший наилучшие результаты при решении задач рассматриваемого класса), но текущая реализация позволяет применять и другие известные варианты метода, в частности, Hestenes--Stiefel, Fletcher--Reeves, Polak--Ribière--Polyak(plus), Con\-ju\-gate Descent--Flet\-cher, Liu--Storey, Dai--Yuan (детали изложены в \cite{andrei_40_CG_2008}).

Достоинством метода, определившего его высокую популярность на практике, является высокая скорость работы и универсальность применения --- при правильной настройке он показывает конкурентоспособные результаты на широком классе задач. К недостаткам стоит отнести необходимость ``рестартов'', т.е. периодическому сбросу истории. К сожалению, на сегодняшний день не существует надёжных критериев рестарта, поэтому на практике применяются различные эвристические подходы, самый простой и популярный из которых --- рестарт каждые $K$ итераций. От того, насколько адекватно настроен критерий рестарта напрямую зависит эффективность работы алгоритма: как при слишком маленьком значении $K$ (тогда метод приближается к методу наискорейшего спуска), так и при слишком большом его значении скорость сходимости падает. Поэтому для получения максимальной производительности на конкретном классе задач следует потратить некоторое время на настройку параметров метода.

Описание реализованного метода CG представлено в алгоритме~\ref{CG_ALGORITHM}.

\begin{algorithm}[!h]
\caption{Метод CG (сопряжённых градиентов)}
\label{CG_ALGORITHM}
\begin{algorithmic}
\REQUIRE $x^0 \in \mathbb{R}^n$
\STATE $r^0 \leftarrow -\nabla f(x^0)$
\STATE $k \leftarrow 0$
\REPEAT
  \STATE $\tilde{r}^k \leftarrow r^k / \| r^k \|_2$
  \STATE $h_k \leftarrow \argmin\limits_{h \ge 0} f( x^k + h \tilde{r}^k )$
  \STATE $x^{k+1} \leftarrow x^k + h_k \tilde{r}^k$
  \STATE $s^{k} \leftarrow x^{k+1} - x^k$
  \STATE $y^{k} \leftarrow \nabla f( x^{k+1} ) - \nabla f( x^k )$
  \IF{\textbf{restart}}
    \STATE \COMMENT{Выполяем рестарт метода (см. описание)}
    \STATE $r^{k+1} \leftarrow -\nabla f(x^{k+1})$
  \ELSE
    \STATE \COMMENT{Вычисляем значение $\beta_k$ согласно используемому варианту метода (см. пункт \ref{CG})}
    \STATE $r^{k+1} \leftarrow -\nabla f( x^{k+1} ) + \beta_k r^k$
  \ENDIF
  \STATE $k = k + 1$
\UNTIL \textbf{convergence}.
\RETURN $x^{k}$
\end{algorithmic}
\end{algorithm}

\subsubsection{LBFGS}
\label{LBFGS_DESC}

Метод LBFGS (Limited-memory BFGS method) является реализацией квазиньютоновского метода BFGS (см. пункт \ref{CG}). Этот метод также является одним из самых популярных на практике в силу своей универсальности и высокой эффективности. Вычислительная сложность его итерации существенно меньше, чем у ``полных'' квазиньютоновских методов и сопоставима с методами FGM и CG, т.к. используются только операции линейной алгебры над векторами, а операции матричного типа отсутствуют. Преимуществом  перед методом CG является отсутствие процедур рестарта, поскольку размер используемой истории явно задаётся параметром $m$. Вычислительный опыт показывает, что в большинстве случаев значения $m$ порядка 3--5 вполне достаточно  и дальнейшее его увеличение лишь повышает вычислительную сложность метода без улучшения качества его работы.

Описание реализованного метода LBFGS представлено в алгоритмах~\ref{LBFGS_DIR_ALGORITHM}~и~\ref{LBFGS_ALGORITHM}.

\begin{algorithm}[!h]
\caption{Вычисление направления спуска (метод LBFGS)}
\label{LBFGS_DIR_ALGORITHM}
\begin{algorithmic}
\REQUIRE $k \ge 0$
\STATE $r \leftarrow -\nabla{f(x^k)}$
\IF{$k = 0$}
  \RETURN $r \leftarrow r / \| r \|_2$
\ENDIF
\FOR{$i = (k-1), (k-2), ..., (k-m)$}
  \STATE $\alpha_i \leftarrow \rho_i \langle s^i, r \rangle$
  \STATE $r \leftarrow r - \alpha_i y^i$
\ENDFOR
\STATE $r \leftarrow \langle s^k, y^k \rangle / \langle y^k, y^k \rangle r$
\FOR{$i = (k-m), (k-m+1), ..., (k-1)$}
  \STATE $\beta \leftarrow \rho_i \langle y^i, r \rangle$
  \STATE $r \leftarrow r + (\alpha_i - \beta) s^i$
\ENDFOR
\RETURN $r$
\end{algorithmic}
\end{algorithm}

\begin{algorithm}[!h]
\caption{Метод LBFGS}
\label{LBFGS_ALGORITHM}
\begin{algorithmic}
\REQUIRE $x^0 \in \mathbb{R}^n$, $m \in \mathbb{N}$
\STATE $k \leftarrow 0$
\REPEAT
  \STATE Вычисляем текущее направление спуска $r^k$ используя алгоритм \ref{LBFGS_DIR_ALGORITHM}
  \STATE $h_k \leftarrow \argmin\limits_{h \ge 0} f\left( x^k + h r^k \right)$
  \STATE $x^{k+1} \leftarrow x^k + h_k r^k$
  \STATE $s^k \leftarrow x^{k+1} - x^k$
  \STATE $y^k \leftarrow \nabla f( x^{k+1} ) - \nabla f( x^k )$
  \STATE $\rho_k \leftarrow 1 / \langle s^k, y^k \rangle$
  \STATE Удаляем вектора $s^{k-m}$, $y^{k-m}$ из памяти
  \STATE $k \leftarrow k + 1$
\UNTIL \textbf{convergence}.
\RETURN $x_k$
\end{algorithmic}
\end{algorithm}

\subsection{Методы одномерного поиска}
\label{LS_METHODS}

Все представленные ранее методы для своей работы используют процедуру одномерного поиска. Рассматриваемый в работе класс оптимизируемых функций является ``вычислительно затратным'', т.е. каждое обращение (запрос) к оракулу (вычисление функции и, особенно, её градиента), требует выполнения достаточно большого объёмы вычислений. Поскольку процесс одномерного поиска неизбежно приводит к таким обращениям и, соответственно, вычислительным затратам, первостепенным вопросом с точки зрения повышения эффективности работы алгоритмов становится необходимость уменьшения количества таких запросов. Использование плохо оптимизированных процедур одномерного поиска, которые многократно вычисляют значение оптимизируемой функции в процессе своей работы, может радикально замедлить общее быстродействие и ``свести на нет'' все усилия, затраченные на вышестоящие алгоритмы. Можно утверждать, что в нашем случае именно процедуры одномерного поиска являются ``бутылочным горлышком'', от которого, в итоге, и зависит скорость поиска локального экстремума.

В рамках проведённой работы авторами реализованы два метода одномерного поиска, представленных в алгоритмах \ref{LS_H_ALGORITHM} и \ref{LS_PAR_ALGORITHM}. Они основаны на хорошо зарекомендовавших себя на практике подходах (см., например \cite{nocedal_wright_2006} \cite{numerical_recipes_2007}), но также реализуют ряд идей, направленных как на учёт специфики исследуемого класса задач, так и на повышение вычислительной эффективности. В частности, используется история работы алгоритма, т.е. шаг одномерного поиска, найденный на предыдущей итерации градиентного метода, становится стартовым значением при следующем запуске одномерной минимизации. Даже если он является ``удачным'', методы одномерного поиска пытаются найти другой шаг, обеспечивающий лучшую релаксацию оптимизируемой функции. Такое поведение страхует нас от ситуаций, когда, например, был (неудачно) выбран небольшой стартовый шаг ($10^{-3}$, например), постоянно применяя который, градиентный метод может спускаться, но скорость такого спуска будет абсолютно неудовлетворительной. В целях повышения быстродействия число обращений к оракулу сокращено до минимально возможного --- в большинстве случаев оказывается достаточно 2-3 вызовов для нахождения релаксирующей точки. Краеугольной идеей предложенных алгоритмов является переход от максимального качества работы (т.е. ``честный'' $\argmin$) к максимальному быстродействию. Проведённые вычислительные эксперименты показали, что за фиксированное время более выгодным оказывается выполнение б\'{о}льшего числа ``грубых'' итераций градиентного метода чем меньшего числа высокоточных.

\begin{algorithm}[!h]
\caption{Метод одномерного поиска LS\_H}
\label{LS_H_ALGORITHM}
\begin{algorithmic}
\REQUIRE $x^0 \in \mathbb{R}^n$, $r \in \mathbb{R}^n$, 
$h_0 > 0$, $0 < \varepsilon_h < 1$, 
$k_{(+)} > 1$, $0 < k_{(-)} < 1$
\STATE $f_* \leftarrow f(x^0)$ \COMMENT{Выполняем, только если значение не вычислено в вызывающем методе}
\STATE $\tilde{f}_0 \leftarrow f(x^0 + h_0 r)$
\IF{$\tilde{f}_0 < f_*$}
  \STATE $h_1 \leftarrow k_{(+)} h_0$
  \STATE $\tilde{f}_1 \leftarrow f(x^0 + h_1 r)$
  \IF{$\tilde{f}_1 < \tilde{f}_0$}
    \RETURN ($h_1, \tilde{f}_1$)
  \ENDIF
  \RETURN ($h_0, \tilde{f}_0$)
\ENDIF
\STATE $h_2 \leftarrow k_{(-)} h_0$
\STATE $\tilde{f}_2 \leftarrow f(x^0 + h_2 r)$
\WHILE{$\tilde{f}_2 > f_*$}
  \STATE $h_2 \leftarrow k_{(-)} h_2$
  \IF{$h_2 \le \varepsilon_h$}
    \RETURN ($0, f_*$) \COMMENT{Невозможно найти релаксирующий шаг}
  \ENDIF
  \STATE $\tilde{f}_2 \leftarrow f(x^0 + h_2 r)$
\ENDWHILE
\RETURN ($h_2, \tilde{f}_2$)
\end{algorithmic}
\end{algorithm}

\begin{algorithm}[!h]
\caption{Метод одномерного поиска LS\_PAR}
\label{LS_PAR_ALGORITHM}
\begin{algorithmic}
\REQUIRE $x^0 \in \mathbb{R}^n$, $r \in \mathbb{R}^n$,
$h_0 > 0$, $K \ge 2$, $G_0 \in \{ true, false\}$
\STATE $\tilde{h}_0 \leftarrow 0$
\STATE $\tilde{f}_0 \leftarrow f_* \leftarrow f(x^0)$ \COMMENT{Значение уже вычислено в вызывающем методе}
\IF{$G_0 = true$}
  \STATE $g_* \leftarrow \nabla f(x^0)$ \COMMENT{Значение уже вычислено в вызывающем методе}
  \STATE $\tilde{h}_1 \leftarrow h_0$
  \STATE $\tilde{f}_1 \leftarrow f(x^0 + \tilde{h}_1 r)$
  \STATE \COMMENT{По 3-м парам $(0, f_*)$, $(0, \langle g_*, r \rangle )$ и $(\tilde{h}_1, \tilde{f}_1)$ строится аппроксимация параболы, точка минимума которой заносится в $\tilde{h}_2$}
  \STATE $\tilde{f}_2 \leftarrow f(x^0 + \tilde{h}_2 r)$
\ELSE
  \STATE $\tilde{h}_1 \leftarrow -1/2 h_0$
  \STATE $\tilde{f}_1 \leftarrow f(x^0 + \tilde{h}_1 r)$
  \STATE $\tilde{h}_2 \leftarrow  1/2 h_0$
  \STATE $\tilde{f}_2 \leftarrow f(x^0 + \tilde{h}_2 r)$
\ENDIF
\FOR{$k = 2, ..., K$}
  \STATE \COMMENT{По 3-м парам $(\tilde{h}_0, \tilde{f}_0)$, $(\tilde{h}_1, \tilde{f}_1)$ и $(\tilde{h}_2, \tilde{f}_2)$ строится аппроксимация параболы, точка минимума которой заносится в $\tilde{h}_3$}
  \STATE $\tilde{f}_3 \leftarrow f(x^0 + \tilde{h}_3 r)$
  \STATE \COMMENT{Вектора $\tilde{h}$ и $\tilde{f}$ согласовано сортируются по возрастанию элементов $\tilde{f}$}
\ENDFOR
\IF{$K = 2$}
  \STATE \COMMENT{Вектора $\tilde{h}$ и $\tilde{f}$ согласовано сортируются по возрастанию элементов $\tilde{f}$}
\ENDIF
\IF{$\tilde{f}_0 < f_*$}
  \RETURN ($\tilde{h}_0, \tilde{f}_0$)
\ENDIF
\RETURN ($0, f_*$) \COMMENT{Релаксирующий шаг не найден}
\end{algorithmic}
\end{algorithm}

\subsection{Вычислительные эксперименты}

Проверка описанных выше методов оптимизации производилась в контексте алгоритмов для решения задачи предсказания белок-белковых комплексов (белок-белковый докинг). Простыми словами данная задача может быть описана следующим образом: для имеющейся пары молекул нужно определить наиболее вероятное взаимное расположение этих молекул при образовании молекулярного комплекса белок-белок. Задача считается успешно решенной, если на выходе имеется структура комплекса, которая соответствует встречаемому в природе комплексу или очень близка к этому положению. Задача белок-белкового докинга имеет большую значимость для решения прикладных задач, связанных с разработкой лекарственных средств, основанных на белковых молекулах, так как терапевтический эффект обусловлен связыванием терапевтического белка с мишенью, а это и есть белок-белковый комплекс. Также стоит отметить, что в рамках типичного проекта по разработке лекарственного препарата, основанного на антителах, в компании БИОКАД эту задачу требуется решать сотни раз.

На промежуточных этапах решения задачи белок-белкового докинга могут быть найдены, и очень часто находятся, физически-невозможные структуры комплексов, в частности, структуры со столкнувшимися атомами. Такие ситуации возникают в следствии использования грубых метрик на начальных этапах с целью ускорения процесса сканирования большого пространства решений. Следствием этой проблемы является не возможность применения наиболее физически обоснованной метрики для финального ранжирования решений -- потенциальной энергии белкового комплекса. Эта метрика является наиболее обоснованной, так как меньшему значению потенциальной энергии соответствует более энергетически выгодная молекулярная структура, отсюда следует что минимуму (либо одному из возможных минимумов) потенциальной энергии соответствует структура, встречаемая в природе. Учитывая вышесказанное, минимизация энергии является одним из важнейших этапов решения задачи докинга, так как позволяет привести структуру комплекса в физически-оправданную конформацию, а значит, и позволяет использовать энергетические метрики в качестве метрик оценки решений, что, в свою очередь существенно повышает точность предсказаний. 

Для оценки методов минимизации был использован тестовый набор для оценки точности белок-белкового докинга. Этот набор состоит из белок-белковых комплексов с достоверно известной трехмерной структурой, встречаемой в природе (нативная структура), другими словами правильный ответ уже известен. Часть данных были получены из прошедших раундов соревнований по белок-белковому докингу CAPRI (Critical Assessment of PRedicted Interactions) \cite{janin2002welcome} и из курируемого набора данных для белок-белкового докинга Protein-Protein Docking Becnchmark 5.0 \cite{vreven2015updates}. Для каждого комплекса с помощью алгоритмов докинга, разработанных в компании БИОКАД, были сгенерированы 499 начальных приближений, отличающихся друг от друга расположением и поворотом в трехмерном пространстве одной из молекул. То есть одна из молекул комплекса оставалась неподвижной, изменялось же расположение только второй молекулы. Кроме того, было искусственно добавлено заведомо верное расположение молекул, соответствующее нативным структурам. Так как нативной структуре комплекса соответствует минимум потенциальной энергии (E), в данном эксперименте производилась минимизация всех начальных приближений, сортировка по значению E и оценка порядкового номера нативной (или очень близкой к таковой) структуры. Если близкая к нативной структура присутствует в первых 30 результатах, задача считалась успешно решенной.

Приводимые ниже результаты численных экспериментов получены на вычислительной системе, имеющей следующие характеристики:

\begin{itemize}
\item Ubuntu server 16.04, x86\_64
\item 2 x Intel Xeon E5-2687Wv3
\item 4 x Nvidia Tesla K80
\item Компилятор gcc 5.4.0 с параметрами:\\
\texttt{--std=c++11 -O2 -Wall -Wextra -Wpedantic -Wdouble-promotion\\ -Wfloat-conversion -s -DNDEBUG}
\item CUDA-компилятор nvcc v9.0.176 с параметрами:\\
\texttt{--std=c++11 -O3 -lineinfo -Xcompiler -Wno-attributes \\-gencode arch=compute\_37,code=sm\_37 -DNDEBUG}
\end{itemize}

В качестве основного режима работы реализованных алгоритмов был выбран float (одинарная точность), что обусловлено как более высоким быстродействием современных GPU при работе с вещественными числами одинарной точности, так и тем, что для требуемых расчетов float-арифметика обеспечивает достаточную точность результатов.

\subsection{Результаты вычислительных экспериментов}

В таблицах \ref{table_res_1} и \ref{table_res_2} представлены результаты вычислительных экспериментов, проведенных на тестовых данных белок-белкового докинга. Аббревиатуры методов соответствуют: LBFGS -- метод LBFGS c параметром $m = 3$ и одномерным поиском из алгоритма \ref{LS_PAR_ALGORITHM}, FGM -- метод FGM c одномерным поиском из алгоритма \ref{LS_H_ALGORITHM}, CG\_PRP -- метод CG в варианте Полака--Рибьера--Поляка c одномерным поиском из алгоритма \ref{LS_PAR_ALGORITHM}. Таблица \ref{table_res_1} описывает наименьший порядковый номер (индекс) структуры, близкой к нативной, при сортировке результатов по значению потенциальной энергии. Чем меньше это число, тем точнее решена задача докинга в данном случае, при этом, если оно меньше 30, то задача считается успешно решенной. Также, приведены средние значения наименьших индексов структур, близких к нативным, по которым можно судить о среднем качестве работы метода на всем тестовом наборе. Для определения близости к нативной структуре была использована метрика RMSD (Root-mean-square deviation), которая является стандартной для измерения расстояния между подобными молекулярными структурами 
\[
RMSD = \sqrt{\frac{1}{n}\sum\limits_{i=1}^n ||v_i - w_i||^2},
\]
где $n$ -- утроенное число атомов (каждый атом описывается тремя координатами), а $v_i$, $w_i$ -- соответствующие координаты соответствующих атомов, сравниваемых структур.
Структура считалась близкой к нативной при значении RMSD < 10Å. Этот порог является максимальным при оценке решений во всемирных соревнованиях по белок-белковому докингу CAPRI \cite{janin2002welcome}, поэтому он и был взят. Значения RMSD также приведены в таблице \ref{table_res_1}. В таблице \ref{table_res_2} представлены значения потенциальной энергии для структуры с наименьшим значением энергии $E_{best}$ и структуры, которая является лучшей из близких к нативной $E_{native}$ (в случае, если это одна и та же структура, приведены одинаковые значения). 

Общее время вычислительного эксперимента для трех методов занимало порядка недели в режиме параллельной работы на 4-х Nvidia Tesla K80 GPU, полная конфигурация сервера описана выше. Общее количество минимизированных структур за один эксперимент составляет 91500: 500 начальных приближений для каждого из 61 комплекса 3-мя различными методами. Реальное время работы алгоритмов особенно важно, поскольку минимизация является наиболее вычислительно массивным этапом белок-белкового докинга, а в типичном процессе разработки лекарственного препарата задачу докинга, как уже упоминалось, нужно решать сотни раз. Текущая реализация позволяет произвести минимизацию 500 начальных приближений менее, чем за час на 4-х Nvidia Tesla K80 GPU, что является приемлемым временем работы.

На основании средних значений наименьших индексов структур, близких к нативным, из таблицы \ref{table_res_1} можно видеть, что методы FGM и CG\_PRP показывают результаты лучше, нежели LBFGS. Тем не менее, тестовый набор данных недостаточно велик для окончательного выбора одного метода. В таблице \ref{table_res_2} можно видеть много случаев, в которых метод LBFGS показывает лучшие результаты с точки зрения абсолютного значения потенциальной энергии, тем не менее, с точки зрения решения практической задачи, решаемой в этом эксперименте, метод LBFGS ощутимо проигрывает.
\begin{table}[h]
  \scriptsize
  \begin{center}
  \begin{tabular}{ l | r | r | r | c | c | c }
    идентификатор 
    & \multicolumn{3}{p{3.5cm}|}{\centering индекс структуры, \\ близкой к нативной} 
    & \multicolumn{3}{p{3.5cm}}{\centering RMSD \\ до нативной позы, Å} \\ 
    \hline
     & LBFGS & FGM & CG\_PRP & LBFGS & FGM & CG\_PRP \\
     capri\_8\_22 & 2 & 0 & 0 & 2.16 & 0.63 & 0.63 \\
     capri\_7\_21 & 4 & 7 & 13 & 0.45 & 1.79 & 1.79 \\
     capri\_5\_17 & 121 & 12 & 20 & 0.82 & 0.87 & 0.87 \\
     capri\_5\_14 & 1 & 0 & 0 & 0.91 & 3.00 & 3.00 \\
     capri\_34\_104 & 1 & 0 & 1 & 0.38 & 1.88 & 1.88 \\
     capri\_32\_98 & 0 & 0 & 0 & 0.38 & 0.88 & 0.88 \\
     capri\_32\_101 & 1 & 0 & 1 & 0.84 & 0.80 & 0.80 \\
     capri\_30\_94 & 2 & 0 & 0 & 0.67 & 0.33 & 0.33 \\
     capri\_30\_88 & 2 & 1 & 0 & 0.35 & 1.20 & 1.20 \\
     capri\_30\_86 & 233 & 173 & 262 & 1.73 & 0.34 & 0.34 \\
     capri\_30\_85 & 0 & 0 & 0 & 0.93 & 0.27 & 0.27 \\
     capri\_30\_84 & 1 & 0 & 0 & 0.29 & 0.86 & 0.86 \\
     capri\_30\_79\_3 & 52 & 311 & 114 & 0.43 & 0.33 & 0.33 \\
     capri\_30\_79\_2 & 8 & 8 & 19 & 0.76 & 0.72 & 0.72 \\
     capri\_30\_79\_1 & 0 & 0 & 0 & 0.92 & 0.35 & 0.35 \\
     capri\_30\_77\_4 & 6 & 2 & 5 & 0.33 & 0.31 & 0.31 \\
     capri\_30\_77\_3 & 211 & 245 & 346 & 0.35 & 0.33 & 0.33 \\
     capri\_30\_77\_2 & 1 & 0 & 0 & 0.35 & 0.34 & 0.34 \\
     capri\_30\_77\_1 & 19 & 6 & 9 & 0.77 & 0.34 & 0.34 \\
     capri\_30\_72\_3 & 14 & 111 & 41 & 0.22 & 0.29 & 0.29 \\
     capri\_30\_72\_2 & 28 & 11 & 98 & 0.23 & 0.29 & 0.29 \\
     capri\_30\_72\_1 & 8 & 7 & 21 & 2.53 & 0.28 & 0.28 \\
     capri\_30\_69\_2 & 233 & 157 & 203 & 0.23 & 0.27 & 0.27 \\
     capri\_30\_69\_1 & 0 & 0 & 0 & 0.29 & 0.30 & 0.30 \\
     capri\_30\_68\_2 & 44 & 23 & 49 & 0.38 & 0.35 & 0.35 \\
     capri\_30\_68\_1 & 4 & 0 & 3 & 0.51 & 0.39 & 0.39 \\
     capri\_27\_58 & 0 & 0 & 1 & 2.21 & 1.04 & 1.04 \\
     capri\_26\_55 & 27 & 37 & 41 & 0.42 & 7.19 & 7.19 \\
     capri\_26\_54 & 316 & 325 & 357 & 0.31 & 0.37 & 0.37 \\
     capri\_24\_50 & 98 & 140 & 35 & 8.10 & 8.05 & 8.05 \\
     capri\_22\_46 & 469 & 485 & 466 & 0.81 & 0.80 & 0.80 \\
     capri\_19\_41 & 27 & 158 & 4 & 1.79 & 0.77 & 0.77 \\
     capri\_18\_40 & 63 & 106 & 317 & 2.35 & 2.35 & 2.35 \\
     capri\_15\_36\_min & 451 & 467 & 436 & 0.40 & 0.49 & 0.49 \\
     capri\_15\_36 & 248 & 409 & 375 & 0.50 & 0.49 & 0.49 \\
     capri\_15\_32\_min & 360 & 29 & 144 & 0.30 & 0.38 & 0.38 \\
     capri\_15\_32 & 13 & 31 & 9 & 3.05 & 3.04 & 3.04 \\
     capri\_11\_27 & 263 & 295 & 183 & 0.47 & 0.41 & 0.41 \\
     capri\_10\_26 & 175 & 81 & 27 & 8.69 & 5.57 & 5.57 \\
     4R9Y & 79 & 67 & 70 & 0.67 & 0.71 & 0.71 \\
     4G6J & 23 & 0 & 1 & 0.33 & 0.38 & 0.38 \\
     4G5Z & 30 & 64 & 87 & 3.78 & 3.75 & 3.75 \\
     3MXW & 5 & 2 & 2 & 0.44 & 0.43 & 0.43 \\
     3EO9 & 188 & 153 & 248 & 0.33 & 0.45 & 0.45 \\
     2JEL & 391 & 143 & 375 & 5.43 & 0.48 & 0.48 \\
     1VFB\_min & 130 & 27 & 3 & 1.76 & 8.77 & 8.77 \\
     1VFB & 49 & 6 & 4 & 1.76 & 8.76 & 8.76 \\
     1PPE & 0 & 8 & 0 & 0.27 & 0.60 & 0.60 \\
     1NTY & 328 & 105 & 55 & 7.61 & 7.65 & 7.65 \\
     1KXQ & 32 & 15 & 25 & 7.62 & 0.41 & 0.41 \\
     1IJK & 165 & 5 & 80 & 0.42 & 0.49 & 0.49 \\
     1I9R & 23 & 0 & 13 & 5.86 & 5.87 & 5.87 \\
     1FSK & 92 & 10 & 93 & 1.53 & 0.50 & 0.50 \\
     1E6J\_min & 178 & 347 & 28 & 1.33 & 1.37 & 1.37 \\
     1E6J & 22 & 32 & 6 & 1.34 & 0.34 & 0.34 \\
     1DQJ & 333 & 262 & 263 & 0.58 & 0.60 & 0.60 \\
     1CGI & 5 & 10 & 0 & 0.62 & 0.62 & 0.62 \\
     1BVK & 101 & 18 & 63 & 0.47 & 9.48 & 9.48 \\
     1BJ1 & 372 & 88 & 238 & 8.76 & 8.84 & 8.84 \\
     1AY7 & 72 & 52 & 77 & 6.77 & 0.52 & 0.52 \\
     1AVX & 5 & 1 & 37 & 0.51 & 0.41 & 0.41 \\
     \hline
     среднее & 100.5 & 82.8 & 88 &
  \end{tabular}
  \end{center}
  \caption{Данные о порядковом номере структуры, близкой к нативной, при ранжировании структур по значению потенциальной энергии белок-белкового комплекса. }
  \label{table_res_1}
\end{table}

\begin{table}[h]
  \scriptsize
  \begin{center}
  \begin{tabular}{ l | r | r | r | r | r | r }
    идентификатор 
    & \multicolumn{2}{c|}{LBFGS}
    & \multicolumn{2}{c|}{FGM} 
    & \multicolumn{2}{c}{CG\_PRP} 
    \\ 
    \hline
     & $E_{best}$ & $E_{native}$
     & $E_{best}$ & $E_{native}$
     & $E_{best}$ & $E_{native}$
     \\
     capri\_8\_22 & -50571.9 & -45926.6 & -45910.7 & -45910.7 & -46031.8 & -46031.8 \\
     capri\_7\_21 & -78483.8 & -78415.2 & -77749.9 & -77649.0 & -78406.9 & -78254.6 \\
     capri\_5\_17 & -92668.5 & -92341.4 & -93201.7 & -93108.6 & -93017.3 & -92909.6 \\
     capri\_5\_14 & -167005.4 & -137625.4 & -137974.8 & -137974.8 & -138476.5 & -138476.5 \\
     capri\_34\_104 & -55082.3 & -52252.9 & -51961.6 & -51961.6 & -52261.8 & -52181.4 \\
     capri\_32\_98 & -91536.9 & -91536.9 & -91371.0 & -91371.0 & -91809.7 & -91809.7 \\
     capri\_32\_101 & -104513.7 & -98930.8 & -98579.6 & -98579.6 & -100455.5 & -99204.5 \\
     capri\_30\_94 & -162684.0 & -160101.2 & -160394.6 & -160394.6 & -160788.5 & -160788.5 \\
     capri\_30\_88 & -79515.5 & -77559.4 & -79207.1 & -77719.3 & -77917.1 & -77917.1 \\
     capri\_30\_86 & -62108.7 & -48289.3 & -48505.1 & -48378.4 & -48692.0 & -48482.8 \\
     capri\_30\_85 & -126991.7 & -126991.7 & -127564.2 & -127564.2 & -128053.6 & -128053.6 \\
     capri\_30\_84 & -118528.3 & -110604.9 & -111089.0 & -111089.0 & -111608.4 & -111608.4 \\
     capri\_30\_79\_3 & -34419.0 & -34325.9 & -34684.9 & -34481.0 & -34799.7 & -34662.9 \\
     capri\_30\_79\_2 & -34484.8 & -34401.5 & -34691.5 & -34675.7 & -34838.9 & -34762.4 \\
     capri\_30\_79\_1 & -34537.4 & -34537.4 & -34852.6 & -34852.6 & -34913.7 & -34913.7 \\
     capri\_30\_77\_4 & -87193.1 & -84520.2 & -84456.4 & -84368.1 & -88233.9 & -84805.8 \\
     capri\_30\_77\_3 & -84664.4 & -84297.0 & -84502.2 & -84140.1 & -93717.8 & -84403.0 \\
     capri\_30\_77\_2 & -86309.6 & -84786.4 & -84639.5 & -84639.5 & -85033.2 & -85033.2 \\
     capri\_30\_77\_1 & -109759.2 & -84593.2 & -84512.3 & -84422.0 & -85030.6 & -84882.7 \\
     capri\_30\_72\_3 & -249637.5 & -202470.5 & -203271.8 & -203002.0 & -203292.8 & -203096.8 \\
     capri\_30\_72\_2 & -210381.7 & -202457.4 & -203257.2 & -203139.8 & -203323.6 & -203021.0 \\
     capri\_30\_72\_1 & -208147.4 & -202503.2 & -203316.9 & -203122.3 & -203224.0 & -203101.1 \\
     capri\_30\_69\_2 & -145319.0 & -140952.8 & -142658.1 & -141713.3 & -142782.6 & -141980.8 \\
     capri\_30\_69\_1 & -141995.4 & -141995.4 & -142747.0 & -142747.0 & -142954.4 & -142954.4 \\
     capri\_30\_68\_2 & -55235.1 & -51297.3 & -51283.9 & -51104.5 & -51750.1 & -51479.0 \\
     capri\_30\_68\_1 & -70131.8 & -51531.3 & -51239.6 & -51239.6 & -74855.2 & -51790.3 \\
     capri\_27\_58 & -72984.6 & -72984.6 & -72805.0 & -72805.0 & -87178.7 & -73048.7 \\
     capri\_26\_55 & -128419.5 & -127281.6 & -139758.6 & -127382.3 & -131501.1 & -127719.4 \\
     capri\_26\_54 & -45554.3 & -45371.6 & -45540.1 & -45374.1 & -45793.4 & -45533.7 \\
     capri\_24\_50 & -128720.1 & -128495.8 & -128904.9 & -128668.2 & -129142.6 & -128981.7 \\
     capri\_22\_46 & -87610.8 & -79320.3 & -87177.7 & -77999.9 & -88032.3 & -79731.6 \\
     capri\_19\_41 & -57530.0 & -57385.3 & -57395.9 & -57235.2 & -57712.4 & -57675.9 \\
     capri\_18\_40 & -76213.3 & -75956.5 & -76568.1 & -76466.1 & -76742.4 & -76493.8 \\
     capri\_15\_36\_min & -112675.1 & -112097.7 & -113226.7 & -112679.8 & -113655.4 & -113161.1 \\
     capri\_15\_36 & -112614.7 & -112311.3 & -113120.4 & -112760.3 & -113613.7 & -113245.4 \\
     capri\_15\_32\_min & -94423.4 & -86603.0 & -86905.2 & -86812.8 & -87088.4 & -86928.4 \\
     capri\_15\_32 & -87701.2 & -86832.0 & -86871.9 & -86759.9 & -87046.8 & -86995.9 \\
     capri\_11\_27 & -78181.6 & -77894.2 & -78157.2 & -77887.1 & -78501.1 & -78248.3 \\
     capri\_10\_26 & -107489.4 & -105406.4 & -105762.5 & -105598.2 & -105937.4 & -105844.0 \\
     4R9Y & -111822.6 & -109022.7 & -108823.6 & -108700.9 & -109864.2 & -109644.2 \\
     4G6J & -119912.2 & -119798.3 & -120160.1 & -120160.1 & -120480.3 & -120470.9 \\
     4G5Z & -156913.6 & -118867.8 & -119250.9 & -119045.6 & -119637.3 & -119498.1 \\
     3MXW & -129724.4 & -126724.2 & -228854.3 & -126604.4 & -127087.7 & -126965.9 \\
     3EO9 & -126967.6 & -126681.8 & -127066.8 & -126838.9 & -127299.4 & -127036.6 \\
     2JEL & -101452.8 & -101051.3 & -101044.3 & -100738.9 & -101482.4 & -101111.0 \\
     1VFB\_min & -74414.2 & -74294.0 & -74494.7 & -74415.6 & -74526.9 & -74517.8 \\
     1VFB & -74986.0 & -74330.7 & -74497.5 & -74441.7 & -74595.2 & -74536.9 \\
     1PPE & -49119.3 & -49119.3 & -49226.0 & -49171.3 & -49256.0 & -49256.0 \\
     1NTY & -110838.1 & -110260.9 & -110457.9 & -110222.1 & -110847.3 & -110669.8 \\
     1KXQ & -143886.8 & -132645.5 & -132951.1 & -132872.6 & -133178.7 & -133115.5 \\
     1IJK & -107683.6 & -107381.5 & -107336.6 & -107287.7 & -109706.4 & -107646.3 \\
     1I9R & -112501.7 & -112357.3 & -112660.1 & -112660.1 & -112942.5 & -112848.2 \\
     1FSK & -121052.2 & -120222.2 & -120037.7 & -119945.2 & -120533.2 & -120344.0 \\
     1E6J\_min & -122041.0 & -98979.1 & -99248.2 & -98994.0 & -99321.4 & -99260.8 \\
     1E6J & -99252.8 & -99068.8 & -99207.9 & -99133.1 & -100978.5 & -99290.5 \\
     1DQJ & -115945.9 & -110563.6 & -111170.0 & -110864.7 & -115192.8 & -111144.9 \\
     1CGI & -62803.7 & -59810.1 & -59746.3 & -59680.5 & -59927.0 & -59927.0 \\
     1BVK & -90899.6 & -72872.2 & -73015.9 & -72801.1 & -80154.7 & -73078.0 \\
     1BJ1 & -109904.1 & -104866.7 & -104420.1 & -104267.4 & -104960.2 & -104675.5 \\
     1AY7 & -40858.6 & -40730.9 & -40723.8 & -40561.6 & -40891.7 & -40739.4 \\
     1AVX & -90394.4 & -81194.0 & -80901.2 & -80891.3 & -81278.9 & -81149.8 \\
  \end{tabular}
  \end{center}
  \caption{Данные о значениях потенциальной энергии белок-белковых комплексов: лучшая энергия из всех начальных приближений ($E_{best}$, кДж/моль) и энергия структуры, близкой к нативной ($E_{native}$, кДж/моль).}
  \label{table_res_2}
\end{table}

Работа А.Ю. Горнова поддержана грантом РФФИ 18-07-00587. Работа А.С. Аникина поддержана грантом РФФИ 18-29-03071 мк. Работа А.В. Гасникова поддержана грантом РНФ 17-11-01027.

\printbibliography

\end{document}